\documentclass[12pt]{article}
\usepackage{graphicx}
\usepackage{amssymb,amsfonts,amsthm}
\usepackage{amsmath,amsthm,amssymb}
\usepackage{amsfonts}

\newtheorem{theorem}{Theorem}[section]
\newtheorem{lemma}[theorem]{Lemma}

\newtheorem{prop}[theorem]{Proposition}
\theoremstyle{definition}

\newtheorem{rk}[theorem]{Remark}

\newcounter{ppp}

\begin{document}
\title{Embedding construction based on  amalgamations of group relators.}
\author{A. Yu. Olshanskii \thanks{The
author was supported in part by the NSF grant DMS 1161294}}
\maketitle

\begin{abstract} An embedding construction $G\hookrightarrow H$ for groups $G$ with length function was
introduced by the author earlier. Here we obtain new properties of this embedding, answering some questions
raised by M.V. Sapir. In particular, an analog of Tits' alternative holds for the subgroups of $H$.
\end{abstract}

{\bf Key words:} group relation, group embedding, free group, van Kampen diagram, small cancellation

\medskip

{\bf AMS Mathematical Subject Classification: } 20F05, 20F06, 20F65, 20E65, 20E07

%\large

\section{ The embedding $\gamma: G\hookrightarrow H$.}

Let ${\cal A}^{\pm 1}=\{a_1^{\pm 1},\dots,a_m^{\pm 1}\}$ be a group alphabet. The length of
arbitrary word $W$ in this alphabet is denoted by $|W|_{\cal A}$ or just by $|W|$.
If the set ${\cal A}$ generate a group $G$ and $W$ is a shortest word in the alphabet ${\cal A}^{\pm 1}$ representing an element $g\in G$,
then $|g|=|g|_{\cal A}=|W|_{\cal A}.$

Consider now a group $G$ with an arbitrary {\it length function}, i.e., with a function $\ell:  G\to {\mathbb N}=\{0,1,\dots\}$
satisfying the conditions

$\bullet$
$\ell(g)=\ell(g^{-1})$ for all $g\in G$, and $\ell(g)=0$ iff $g=1$;

$\bullet$     $\ell(gh)\le \ell(g)+\ell(h)$ for $g,h\in G$;

$\bullet$ there exists a positive number $c$ such that
      $\,card\{g\in G\mid \ell(g)\le r\}\le c^r$ for any $r\in {\mathbb N}$.

Obviously, $G$ is a left-invariant metric space with respect to $dist(g,h)=\ell(g^{-1}h)$,
and it is clear that the length with respect to a finite set of generators satisfies
these three conditions. Moreover, if $G$ is a subgroup of a finitely generated
group $H=\langle\cal A\rangle$, then the restriction of the length $|\;|_{\cal A}$
to $G$ is a length function on $G$. Up to equivalence such length functions
do not depend on the choice of a finite generating set of $H$.
Two functions $f_1, f_2: G\to {\mathbb N}=\{0,1,\dots\}$ on a group $G$ are called
{\it equivalent} here if there is a positive constant $c$ such that $f_1(g)\le cf_2(g)$
and $f_2(g)\le cf_1(g)$ for every $g\in G.$

One can show that there are uncountably many pairwise nonequivalent lengths functions
on every infinite countable group (see Corollary 1 in \cite{BO}). Nevertheless, it was proved
that every length function on a group $G$ can be obtained by means of an embedding into a finiely
generated group:

\begin{prop} \label{99} (\cite{O99}) Assume that $\ell$ is a length function on a group $G$. Then there is an
embedding of $G$ into a group $H$ with two generators $a$ and $b$, such the restriction of
the function $g\to |g|_{\{a,b\}}$  to the subgroup $G$ is equivalent to the function $\ell.$
\end{prop} $\Box$

This embedding theorem helped to answer a number of questions raised earlier (see \cite{O99}, \cite{O97}, \cite{OS}).
Here we establish some additional properties of the construction. These characteristics are useful in connection
with ``rapid decay" (RD) property of groups, and the author have been asked about their validity by Mark Sapir.
The results of the present paper are needed in \cite{S} to obtain $2$-generated groups without RD, where all amenable
subgroups are cyclic and undistorted.

At first we want to describe the embedding from \cite{O95} and \cite{O99}.
The notation $U\equiv V$ means that the words $U$ and $V$ in a group alphabet ${\cal A}^{\pm 1}$ are letter-by-letter
equal; these two words are {\it freely equal} ({\it freely conjugate}) if they represent the same
element (the conjugate elements) of the free group $F({\cal A})$.

We shall say that a subset of words
${\cal X}$ is an {\it exponential} set if there exist constants
$C$ and $c>1$ such that $\,card\{X\in{\cal X}|\,|X|\le i\}\ge c^i$ for
every $i\ge C$.

A word $W$ in the `positive' alphabet ${\cal A}=\{a_1,\dots,a_m\}$
is called {\it positive}.

A reduced word $X$ is called $s$-{\it aperiodic} if it has no
non-empty powers $Y^s$ as subwords.

\begin{lemma} \label{exp} For any $\lambda >0$ one can choose an exponential
set $\cal Y$ of positive words over the alphabet
$\{a,b\}$ with the following properties.

(*) Let $V$ be a subword of some $W\in {\cal Y}$ and $|V|\ge \lambda |W|$.
Then $V$ occurs in $W$, as a subword, only once. If this $V$ is a
subword of some $U\in {\cal Y}$, then $U\equiv W$.

(**) Every word in $\cal Y$ is $7$-aperiodic.

\end{lemma}

\proof We recall a construction from \cite{O99}, Lemma 4. Given a $2$-letter alphabet
$\{a,b\}$, there is an exponential set of $6$-aperiodic positive words ${\cal X}_1$
such that every word $X\in {\cal X}_1$ starts and ends with the letter $b$.
We enumerate the word in ${\cal X}_1$ so that $i<j$ if $|X_i|<|X_j|$. It was shown
in Lemma 4 of \cite{O99} that there is an integer $N_0$ such that the set $\cal Y$
consisting of all words
$$Y_j\equiv a_2^6X_{(j-1)N_0+1}a_2^6X_{(j-1)N_0+2}\dots a_2^6X_{jN_0},\;\; j=1,2,\dots, $$
is exponential and satisfies (*). It remains to verify (**).

Assume that $Y_j$ has a nonempty subword $Z$ of the form $A^{7}$. Since every $X_k$
is $6$-aperiodic and starts/ends with $b$, the subwords $a^6X_ka^6$ are $7$-aperiodic.
It follows that $A^{7}$, and so some cyclic shift $A'$  of $A$, must contain the subword
$a^6$. We have a $6$-th power  $(A')^6\equiv (a^6U)^6$ as a subword of $Z$ and of $Y_j$.
But this is not possible since all the words $X_k$ are pairwise distinct, a contradiction.
\endproof

Let us fix $\lambda =0.01$ and let $\cal Y$ be an exponential set of words in the
alphabet $\{a,b\}$ satisfying the properties (*) and (**) from  Lemma \ref{exp}.
Given a group $G$ with a length function $\ell$,
there
are a constant $d=d({\cal Y},\ell)$ and a subset
${\cal X}=\{X_g\}_{g\in G}\subset {\cal Y}$ such that the inequalities
$$\ell(g)\le |X_g|<d\,\ell(g),\,\,\,g\in G\backslash \{1\}$$
hold (see \cite{O99}, Lemma 5; and note that the condition (*) for
some $\lambda$ implies the same condition for any $\mu>\lambda$).

The group $G$ can be presented as a homomorphic image of the free
group $F_G$ with the basis $\{x_g\}_{g\in G\backslash\{1\}}$ under the
epimorphism $\varepsilon :\,x_g\mapsto g$. Denote by $\delta$ the
natural embedding of the kernel $\ker \varepsilon = N$ into the group $F_G$.

Below we explain  the commutative diagram

{
$$ \begin{array}{ccccccccc}
1 & \longrightarrow & N & \stackrel{\mbox{\large$\delta$}}\longrightarrow & F_G &
\stackrel{\mbox{\large$\varepsilon$}}\longrightarrow &G & \longrightarrow & 1 \vspace{0.9em}\\
  &  & \downarrow\alpha &  & \downarrow\beta & & \downarrow\gamma & & \vspace{0.9em}\\
1 & \longrightarrow & L & \stackrel{\mbox{\large$\bar\delta$}}\longrightarrow & F(a,b) &
\stackrel{\mbox{\large$\bar\varepsilon$}}\longrightarrow & H & \longrightarrow & 1 \\
\end{array}
%\eqno{(7)}
$$
} \medskip

Namely, the group $F(a,b)$ is free with the basis  $\{a.b\}$, and the homomorphism $\beta$ is given by the formula
$\beta (x_g)=X_g$ for $g\in G\backslash\{1\}$. By definition, the homomorphism $\bar \delta$
is the canonical embedding of the normal closure $L$ of the subgroup
$\beta \delta (N)$ in $F(a,b)$. The homomorphism $\alpha$ is well-defined
by the equality $\bar\delta\alpha=\beta\delta$. Define $\bar\varepsilon$ as the natural epimorphism
of the group $F(a,b)$ onto its quotient $H=F(a,b)/L$, and so one can
regard the set $\{a,b\}$ as a generating set
for the group $H$ too.
Finally, in view
of the condition $\bar\delta\alpha=\beta\delta$, the homomorphism $\gamma$ can
be well-defined by the equality $\gamma\varepsilon=\bar\varepsilon\beta$.

By Lemma 6 \cite{O99}, the homomorphism $\gamma$ is injective, and so one can identify
the group $G$ with the subgroup $\gamma(G)\le H$, and by Theorem 1 \cite{O99}, the
function  $g\mapsto |g|_{\{a,b\}}$ is equivalent to the function $\ell$.
(Only condition (*) with $\lambda=0.02$ was used for $\cal X$ in \cite{O99} to obtain the embedding $\gamma$.)
Here we show that the subgroups of $H$ 'avoiding' the conjugates of $G$ are big, with
few trivial exceptions.

\begin{theorem} \label{t1} Every subgroup $K$ of $H$ is either

(1) conjugate to a subgroup of $G$, or

(2) infinite cyclic and has trivial intersection with any  $hGh^{-1}$ ($h\in H$) , or

(3) infinite dihedral (the involutions from $H$ are conjugate to the elements of $G$), or

(4) contains a noncyclic free subgroup.
\end{theorem}

Every element, which is not conjugate with an element of $G$, generates
a quasi-geodesic cyclic subgroup:

\begin{theorem} \label{t2} If an element $g$ of $H$ is not conjugate with an element of the
subgroup $G$, then $|g^n|_{\{a,b\}}>cn|g|_{\{a,b\}}$ for some $c=c(g)>0$ and every $n>0$.
\end{theorem}

The following theorem implies that the infinite cyclic subgroup $\langle g\rangle $ from Theorem \ref{t2}
 has the Morse property. i.e., the obstacles of linear size intercepting a path along this subgroup cannot be circumvented in linear time. So any asymptotic cone of the group $H$ has cut points \cite{DMS}.

 \begin{theorem} \label{t4}  If an element $g$ of $H$ is not conjugate with an element of the
subgroup $G$ and $\theta\in (0, 1)$, then the lengths of the paths $q$ connecting the
vertices $g^n$ and $g^{-n}$ in the Cayley graph of $H$ and avoiding the ball of radius
$\theta n$ centered at $1$, are not bounded from above by a linear function of $n$.
Moreover, there is $c=c(g)>0$ such that $|q|>cn^2$ for every such a path $q$.
\end{theorem}

The next theorem implies both the malnormality of the embedding $\gamma$ (take $n=2$ in part (2))
and the congruence extension property of the embedded subgroup $G$ :
Every normal in $G$ subgroup $A$ is an intersection $G\cap B$, where the subgroup $B$
is normal in $H$ (choose $B$ as the normal closure of $A$ in $H$ and apply part (1)).
To formulate the theorem, we say that a family of elements $(g_1,\dots,g_n)$ is {\it reluctant}
in a group $G$ if there is no subfamily $(g_{i_1},\dots,g_{i_m})$, where $m\in [1,\dots,n-1]$,
such that the product of conjugacy classes $\prod_{j=1}^m g_{i_j}^G$ contains $1$.

\begin{theorem} \label{t3}  (1) If $g_1,\dots, g_n\in G$ and the equation $\prod_{i=1}^n x_ig_ix_i^{-1}=1$
has a solution in $H$, then it has a solution in $G$.

(2)  If $n\ge 1$ and for a reluctant in $G$ family $(g_1,\dots,g_n)$,  we have $\prod_{i=1}^n x_ig_ix_i^{-1}=1$
in $H$,
then all the elements $x_1,\dots, x_n$ belong to a single left coset of $G$ in $H$.
 \end{theorem}

\section{Amalgamations of relators, small cancellations and diagrams.}

We consider the following set $\cal R$ of defining words for the
presentation $H=F(a,b)/L$.
For any nonempty cyclically reduced word $R=R(x_g,\dots,x_h)$ over the alphabet
$\{x_g\}_{g\in G\backslash\{1\}}$, vanishing in $G$, we include the word $R(X_g,\dots,X_h)$ over the alphabet $\{a^{\pm 1},b^{\pm 1}\}$
into the set $\cal R$. It is clear from the definition of the subgroup $L$ that
$L$ is the normal closure of the set $\cal R$ in the group $F(a,b).$

It follows from the
condition (*) that the set $\cal X$ is a free basis of the subgroup $\beta(F_G)$ (see Lemma 2 in \cite{O95}).
In general, any word in the alphabet $\{a^{\pm 1}, b^{\pm 1}\}$ of the form
\begin{equation}\label{ent}
A_{1}\dots A_{k},
\end{equation}
where $A_i\equiv X^{\pm 1}_{g_i}$ for some $g_i\in G$, will be called a $G$-word. A $G$-word will be considered
with its decomposition (\ref{ent}) into the factors $A_1,\dots, A_k$, which we call {\it entire factorization}.

An $\cal R$-word is, by definition, a cyclically reduced form (over the alphabet
$\{a^{\pm 1}, b^{\pm 1}\}$) of a word from $\cal R$. Note that possible cancellations between neighbor entire factors
$A_i$ and $A_{i+1}$ are small by the condition (*), provided $A_{i+1}\ne A_i^{-1}.$

Below, $G$-words will be viewed as certain paths in the Cayley graph of $H$ or in a van Kampen diagram.
Suppose $p=e_1\dots e_n$ is such a path of length $|p|=n$, the edges $e_s$ of which are labelled by
the symbols $Lab(e_s)\in \{a^{\pm 1},b^{\pm 1}\}$. If a $G$-word happens to be the label $Lab(p)$,
then $p$ splits into a product $p=p_1\dots p_k$ with labels $Lab(p_i)\equiv A_i\equiv X^{\pm 1}_{g_i}$.
In this case, the factorization $p_1\dots p_k$ will be called the {\it entire factorization} of $p$,
and its vertices, which divide $p$ into the {\it segments} $p_i$, will be called {\it entire vertices} of $p$.

We use van Kampen's interpretation (see \cite{LS},\cite{O89}) of the deduction of consequences from defining
relations, according to which for any word $w=w(a,b)$ in the normal closure $L$ of $\cal R$ in $F(a,b)$,
there exists a finite connected and simply-connected, planar $2$-complex (= {\it disc diagram}) $\Delta$, the label of
each edge of which is a letter from $\{a^{\pm 1}, b^{\pm 1}\}$, the label of the boundary contour $\partial\Pi$ of each
2-cell (or just {\it cell}) $\Pi$ is a word from $\cal R$, and the word $w$ is written on the outer contour $\partial\Delta$.

However, successive edges $e,f$ that are incident with an entire vertex of the contour $\partial\Pi$ can have mutual
inverse labels, which enables us to assume that they are mutual inverse in the contour $\partial\Pi.$ By
condition (*), under such an identification one can contract at the end of each segment in $\partial\Pi$,
at most $0.01$ of it. The cyclically reduced path with label $B_1\dots B_k$ (where
$B_i$ is a subword of $A_i$) will be called the {\it reduced contour} $\partial'\Pi$ of $\Pi.$ So $\partial'\Pi$
bounds a region containing the cell $\Pi$ and all the edges from $\partial\Pi\backslash \partial'\Pi$ as well (Fig. 1).

\begin{figure}[h!]
%\label{fig1}
\begin{center}
% This is a LaTeX picture output by TeXCAD.
% File name: [emb1.pic].
% Version of TeXCAD: 4.3
% Reference / build: 30-Jun-2012 (rev. 105)
% For new versions, check: http://texcad.sf.net/
% Options on the following lines.
%\grade{\on}
%\emlines{\off}
%\epic{\off}
%\beziermacro{\on}
%\reduce{\on}
%\snapping{\off}
%\pvinsert{% Your \input, \def, etc. here}
%\quality{8.000}
%\graddiff{0.005}
%\snapasp{1}
%\zoom{4.0000}
\unitlength 1mm % = 2.845pt
\linethickness{0.4pt}
\ifx\plotpoint\undefined\newsavebox{\plotpoint}\fi % GNUPLOT compatibility
\begin{picture}(166,65.25)(12,30)
\thicklines
%\emline(20.75,72.25)(21,72.5)
\multiput(20.75,72.25)(.03125,.03125){8}{\line(0,1){.03125}}
%\end
%\emline(78.25,72.25)(78.5,72.5)
\multiput(78.25,72.25)(.03125,.03125){8}{\line(0,1){.03125}}
%\end
%\emline(21,72.5)(39.25,85.25)
\multiput(21,72.5)(.0482804233,.0337301587){378}{\line(1,0){.0482804233}}
%\end
%\emline(90.75,72.5)(109,85.25)
\multiput(90.75,72.5)(.0482804233,.0337301587){378}{\line(1,0){.0482804233}}
%\end
%\emline(137.25,72.5)(155.5,85.25)
\multiput(137.25,72.5)(.0482804233,.0337301587){378}{\line(1,0){.0482804233}}
%\end
%\emline(39.25,85.25)(56.5,73)
\multiput(39.25,85.25)(.0473901099,-.0336538462){364}{\line(1,0){.0473901099}}
%\end
%\emline(109,85.25)(126.25,73)
\multiput(109,85.25)(.0473901099,-.0336538462){364}{\line(1,0){.0473901099}}
%\end
%\emline(56.5,73)(49.5,53.5)
\multiput(56.5,73)(-.033653846,-.09375){208}{\line(0,-1){.09375}}
%\end
%\emline(126.25,73)(119.25,53.5)
\multiput(126.25,73)(-.033653846,-.09375){208}{\line(0,-1){.09375}}
%\end
%\emline(21.25,72.5)(29.5,54)
\multiput(21.25,72.5)(.0336734694,-.0755102041){245}{\line(0,-1){.0755102041}}
%\end
%\emline(91,72.5)(99.25,54)
\multiput(91,72.5)(.0336734694,-.0755102041){245}{\line(0,-1){.0755102041}}
%\end
%\emline(137.5,72.5)(145.75,54)
\multiput(137.5,72.5)(.0336734694,-.0755102041){245}{\line(0,-1){.0755102041}}
%\end
\put(29.5,54){\line(1,0){20.25}}
\put(99.25,54){\line(1,0){20.25}}
\put(145.75,54){\line(1,0){20.25}}
%\emline(29.75,54.25)(32.5,57.5)
\multiput(29.75,54.25)(.03353659,.03963415){82}{\line(0,1){.03963415}}
%\end
%\emline(99.5,54.25)(102.25,57.5)
\multiput(99.5,54.25)(.03353659,.03963415){82}{\line(0,1){.03963415}}
%\end
%\emline(21.5,72.25)(26,71.25)
\multiput(21.5,72.25)(.15,-.0333333){30}{\line(1,0){.15}}
%\end
%\emline(91.25,72.25)(95.75,71.25)
\multiput(91.25,72.25)(.15,-.0333333){30}{\line(1,0){.15}}
%\end
\put(39.25,85.25){\line(0,-1){3.5}}
\put(109,85.25){\line(0,-1){3.5}}
%\emline(47.75,56.5)(49.75,54.25)
\multiput(47.75,56.5)(.03333333,-.0375){60}{\line(0,-1){.0375}}
%\end
%\emline(117.5,56.5)(119.5,54.25)
\multiput(117.5,56.5)(.03333333,-.0375){60}{\line(0,-1){.0375}}
%\end
\put(2.5,60.75){\line(0,-1){.25}}
\put(60,60.75){\line(0,-1){.25}}
\put(28,64.75){\makebox(0,0)[cc]{$A_1$}}
\put(20.25,61.5){\makebox(0,0)[cc]{$B_1$}}
\put(37.25,57){\makebox(0,0)[cc]{$A_2$}}
\put(35.75,51){\makebox(0,0)[cc]{$B_2$}}
\put(31.5,75.75){\makebox(0,0)[cc]{$A_k$}}
\put(25.75,81.5){\makebox(0,0)[cc]{$B_k$}}
\put(44.25,69){\makebox(0,0)[cc]{$\Pi$}}
%\emline(122,72)(131,73.25)
\multiput(122,72)(.23684211,.03289474){38}{\line(1,0){.23684211}}
%\end
%\emline(131,73.25)(142,72)
\multiput(131,73.25)(.28947368,-.03289474){38}{\line(1,0){.28947368}}
%\end
\put(29.25,20.75){\line(1,0){.25}}
%\emline(164.75,55.5)(165,54.75)
\multiput(164.75,55.5)(.03125,-.09375){8}{\line(0,-1){.09375}}
%\end
\put(155.25,85){\line(1,-3){10.25}}
%\emline(163.25,57)(165.25,54.5)
\multiput(163.25,57)(.03333333,-.04166667){60}{\line(0,-1){.04166667}}
%\end
\put(107.75,66.25){\makebox(0,0)[cc]{$\Pi_1$}}
\put(149.5,62){\makebox(0,0)[cc]{$\Pi_2$}}
\put(119.25,71){\makebox(0,0)[cc]{$o_1$}}
\put(143.75,71.5){\makebox(0,0)[cc]{$o_2$}}
%\emline(102.25,81.75)(105.25,83)
\multiput(102.25,81.75)(.07894737,.03289474){38}{\line(1,0){.07894737}}
%\end
%\emline(103.5,80.5)(105.75,82.75)
\multiput(103.5,80.5)(.03358209,.03358209){67}{\line(0,1){.03358209}}
%\end
%\emline(157.75,74)(159.25,71.75)
\multiput(157.75,74)(.03333333,-.05){45}{\line(0,-1){.05}}
%\end
\put(159.75,74.5){\line(0,-1){2.75}}
%\emline(129.25,74)(132,73.5)
\multiput(129.25,74)(.1833333,-.0333333){15}{\line(1,0){.1833333}}
%\end
%\emline(130,72.5)(132.5,73)
\multiput(130,72.5)(.1666667,.0333333){15}{\line(1,0){.1666667}}
%\end
\put(95.5,80.75){\makebox(0,0)[cc]{$q_1$}}
\put(163.5,67.75){\makebox(0,0)[cc]{$q_2$}}
\put(130,70){\makebox(0,0)[cc]{$x$}}
\put(38,40.5){\makebox(0,0)[cc]{Figure 1}}
\put(123.25,41.25){\makebox(0,0)[cc]{Figure 2}}
\end{picture}
\end{center}
\end{figure}

By a $G$-fragment of the (cyclic) word (\ref{ent}), we mean its subword of the form $A'A_s\dots A_{s+t}A''$,
where $A'$ is a suffix of $A_{s-1}$, $A''$ is a prefix of  $A_{s+t+1}$. This factorization is considered
entire.  Analogously, we define a $G$-fragment of an $G$-path with label (\ref{ent}). If  the label of
the boundary of a diagram $\Delta$ (or of a subpath $p$ of $\partial\Delta$) is an $G$-word, then
we also can distinguish entire vertices on $\partial\Delta$ (or on $p$) and introduce the reduced boundary
$\partial'\Delta$ (or $p'$). The parts of $\partial\Delta$ (of $p$) vanishing when passing to $\partial'\Delta$
(to $p'$, resp.) are situated on the plane outside of the region  bounded by $\partial'\Delta.$

As in \cite{O95}, here we use
a modified small cancellation property related to an amalgamation
of cells in van Kampen diagrams rather than to their annihilation. Two different
cells $\Pi_1$ and $\Pi_2$ of a diagram $\Delta$ are {\it compatible}, if there  are
entire vertices $o_1$ on $\partial\Pi_1$, $o_2$ on $\partial\Pi_2$, and a simple path $x=o_1-o_2$
connecting them, whose label is equal in $F(a,b)$ to an $G$-word (Fig. 2). It is easy to
see that this property does not depend on the choice of the pair $(o_1,o_2)$.
Similarly we define the compatibility of a cell with (a fragment of ) $\partial\Delta.$

If the clockwise boundary contours $q_1$ and $q_2$ of the compatible cells $\Pi_1$ and $\Pi_2$ start with
$o_1$ and $o_2$, resp., then the path $q_1xq_2x^{-1}$ bounds a  subdiagram $\Gamma$
with boundary label freely conjugate to a relator from $\cal R$ (or to $1$). So one can replace this
subdiagram by a single cell  (or remove both $\Pi_1$ and $\Pi_2$). Similarly, if $\Pi$ is compatible with a subpath $p$
of $\partial\Delta$ labeled by a $G$-word, one can remove $\Pi$ from $\Delta$,
replacing $p$ by another path labeled by a $G$-word.

A diagram having no pairs of compatible cells is called {\it reduced}.

We need the following property of reduced  diagrams over the
%presentation  ${\cal P}=\langle a,b |\cal R\rangle$ of the
group $H$.

 \begin{lemma} \label{dia} Let $\Delta$ be a reduced diagram over the
presentation $\langle a,b |\cal R\rangle$. Then

(1) there are no two different cells $\Pi_1, \Pi_2$ in $\Delta$ such that $|q|\ge 0.05|\partial'\Pi_1|$
for some common subpath $q$ of the contours $\partial'\Pi_1, \partial'\Pi_2$;

(2) if $q$ is a common subpath for a boundary path $\partial'\Pi =qq_0$ of a
cell $\Pi$ and for a subpath $p=p_1qp_2$ of the contour $\partial'\Delta$, such that
the label $P$ of $p$ is an $G$-word and $|q|\ge 0.05 |\partial'\Pi|$,
then the label of the path $\bar p = p_1(q_0)^{-1}p_2$ is freely equal to a $G$-word
$\bar P$ with $\bar\varepsilon(\bar P)=\bar\varepsilon(P)$.

\end{lemma}

\proof Under the assumption $|q|\ge 0.05|\partial\Pi_1|$, the cells $\Pi_1$ and $\Pi_2$ become
compatible due the condition (*). See details in Lemma 3 of \cite{O95}, but now plug $0.01$
for $\lambda$. The second statement is also explained in that lemma.
\endproof

We say that a (non-oriented) edge $e$ is {\it inner} if it separates two cells in a diagram.
An edge of the boundary $\partial'\Pi$ of a cell $\Pi$  is {\it outer}  if it lies on $\partial'\Delta$.

We also say that $q'$ is an outer arc of $\partial'\Pi$ in a disc diagram $\Delta$ if $q'$
also belongs to the reduced contour $\partial'\Delta$. By Lemma \ref{dia}, a reduced
diagram satisfies the condition $C'(\lambda)$ with $\lambda =0.05$ (see \cite{LS}, Chapter 5).
Therefore by Grindlinger Lemma \cite{G} with $1-3\lambda=0.85$, we have

\begin{lemma}\label{Grin} If a reduced disc diagram $\Delta$ has at least one cell, then

(1) it has a cell $\Pi$ with an outer arc of length $>0.85|\partial'\Pi|$;

(2) the number of inner edges in $\Delta$ is less than $0.15\Sigma,$ where $\Sigma$ is the sum of the reduced perimeters
$|\partial'\Pi|$  of all cells in $\Delta$.
 \end{lemma}
\endproof
We say that the cell and the arc from Lemma \ref{Grin} (1) are {\it Grindlinger} cell and arc.

 \begin{lemma} \label{obe} Let $\Delta$ be a reduced
diagram with contour $p_1q_1p_2q_2$, where the paths $q_1$ and $q_2$ are reduced and
$|p_1|+|p_2|\le 0.01(|q_1|+|q_2|)$. Assume that no edge belongs to both $q_1$ and $q_2^{-1}$
and no cell
$\Pi$ of $\Delta$ has at least $0.55|\partial'\Pi|$ edges belonging to one the paths $q_i,$ $i=1,2.$
Then there is a cell $\Pi$ whose boundary $\partial'\Pi$ has  edges of both $q_1$ and $q_2$,
and the number of edges of $\partial'\Pi$ belonging either to $q_1$ or to $q_2$ is greater than  $0.55|\partial'\Pi|$.
 \end{lemma}
 \proof There are no edges $e$ such that both
$e$ and $e^{-1}$ belong to $q_i$, $i=1,2$. Indeed, otherwise we would have a subdiagram $\Gamma$ bounded by a part
of $q_i$, contrary to Lemma \ref{Grin} (1) applied to $\Gamma,$
because $0.85>0.55$. The inequality $|p_1|+|p_2|\le 0.01(|q_1|+|q_2|)$
implies that at most $0.01(|q_1|+|q_2|)$ non-oriented edges of $q_1$ and $q_2$ are
shared with $p_1$ or $p_2$. Therefore the number of the non-oriented outer edges $Q$ from the reduced boundaries of the cells in $\Delta$,
lying on $q_1$ or on $q_2$, is at least $0.99(|q_1|+|q_2|).$

It follows from Lemma \ref{Grin} (2) that $I<\frac37 O$, where $I$ and $O$ are the numbers of inner and outer edges in the reduced
contours of the cells from $\Delta$, resp.
Also $O\le Q+|p_1|+|p_2|\le \frac{100}{99}Q$, and so $I+|p_1|+|p_2| \le \frac37 \frac{100}{99} Q+ \frac{1}{99} Q<0.45Q.$
Therefore there exists a cell $\Pi$ having $>0.55|\partial'\Pi|$ edges on $q_1$ and $q_2$.
By the assumption, these edges cannot belong only to $q_1$ or only to $q_2$, which proves the lemma.
\endproof

\section{\bf Periodic words with minimal periods.} We say that an element of the group $H$ (or a word $A$ representing
this element) is {free} if it is not conjugate in $H$ to the element of the subgroup $G$. A free
word $A$ is called {\it minimal} if it is not conjugate in $H$ to a shorter word. A word $W$ is
called $A$-periodic (or periodic with period $A$) if it is a subword of some power $A^t$ with $t>0$.

\begin{lemma} \label{1.1} Assume that $A$ is a minimal word and $W$ is an $A$-periodic word.
Let also $W$  be a subword
of a cyclically reduced form $V$ of a $G$-word.

(1) If $|W|\ge 0.55|V|$, then $|W|<1.1|A|$.

(2) We have $|W|<16|A|.$

\end{lemma}

\proof (1) Let $W\equiv B'B_s\dots B_{s+t}B''$, where $B_i$-s are obtained
after (small) cancellations from the factors $A_i$-s of a $G$-fragment. Denote also $B_{s-1}\equiv B'$ and $B_{s+t+1}\equiv B''$.
If $|W|\ge 1.1|A|$, then $W$ can be presented as $UW_1\equiv W_1U'$, where $U$ and $U'$ are
cyclic permutations of $A$ and $|W_1|\ge \frac{1}{11}|W|\ge 0.05|V|$. Now it follows
from (*) that $W_1 \equiv W_2W_3W_4$, where for $i\ne j$,  $W_3$ is a subword
of two factors $B_i$ and $B_j$ with $|W_3|>|W_1|/4>0.011|B_i|$. Therefore $B_i\equiv B_j$
and $W_3$ occurs in $B_i$ only once. Now the equality $UW_2W_3W_4\equiv W_2W_3W_4U'$ implies
that $U$, and so $A$, is freely conjugate to  the $G$-word $A_i\dots A_{j-1}$.
This contradiction with the assumption that $A$ is free, proves the first statement.

(2) Again, the word $W$ has the form $C'C_l\dots C_m C''$, where every factor is a subword of a word from ${\cal X}^{\pm 1}.$
Therefore $|C'|, |C''|<7|A|$ by the property (**). But the middle part $C_l\dots C_m$ is a reduced form of a $G$-word
$V$ with $|C_l\dots C_m|>0.98|V|$ by the property (*). Hence $|C_l\dots C_m|< 1.1|A|$  by the part (1).
Therefore $|W|<(7+7+1.1)|A|$, and the statement (2) is proved.
\endproof

\begin{lemma}\label{0.55} Let $A$ be a minimal word,
$p$ be an outer arc of a cell $\Pi$ in a diagram $\Delta$, and $Lab(p)$
be an $A$-periodic word. Then $|p|<0.55|\partial'\Pi|$.
\end{lemma}

\proof If $|p|\ge 0.55|\partial'\Pi|$, then $|p|< 1.1|A|$  by Lemma \ref{1.1}. Then $p$ can be presented as $p_1p_2$, where
$|p_1|\le |A|$, and so $p_1$ is geodesic, and $|p_2|< 0.1|p|$.
If $pq$ is the contour $\partial'\Pi$, then $|p_1|\le |p_2|+|q|$,
whence $|q|\ge 0.9|p|$. Hence $|q|> 0.45(|q|+|p|)=0.45|\partial'\Pi|$, and so $|p|<(1-0.45)|\partial'\Pi|$, a contradiction.

\endproof

\begin{lemma} \label{1touch} Let $\Delta$ be a reduced diagram.

(1) If $\Gamma$ is a  subdiagram with contour $pq$, where $p$ and $q$ are subpaths
of the reduced boundaries of two cells $\Pi_1$ and $\Pi_2$ of $\Delta$, and $\Gamma$
includes neither $\Pi_1$ nor $\Pi_2$, then $\Gamma$ has no cells, i.e., $p=q^{-1}$.

(2) Let $p$ be a subpath with an $A$-periodic label in $\partial'\Delta$,
where $A$ is a minimal word. If the reduced boundary $\partial'\Pi$
of a cell $\Pi$ from $\Delta$ has two different vertices $o_1$ and $o_2$ from $p$, then the whole
subpath $o_1-o_2$ of $p$ belongs to $\partial'\Pi$.
\end{lemma}
\proof (1) Assume that $\Gamma$ has a cell, and $\Gamma$ is a counter-example with
minimal number of cells. Then $\Gamma$ has a Grindlinger cell $\pi$  by Lemma \ref{Grin} (1). A unique maximal common subpath of $\partial'\pi$ with $p$ (with $q$) can be of length at most
$0.05|\partial'\pi|$ by Lemma \ref{dia} (1).  We obtain a contradiction since $0.05+0.05<0.85.$

 (2) We consider the subdiagram $\Gamma$ situated between $\Pi$ and $p$. It suffices to prove
that $\Gamma$ has no cells. Assume that $\Gamma$ has a cell, and $\Gamma$ is a counter-example with
minimal number of cells. Then $\Gamma$ has a Grindlinger cell $\pi$  by Lemma \ref{Grin} (1). A unique maximal common subpath of $\partial'\pi$ with $p$ has length at most
$0.55|\partial'\pi|$ by Lemma \ref{0.55}. By statement (1),  $x$ has at most one maximal common arc with $\partial'\Pi$,
and its length is $< 0.05|\partial'\pi|$ by Lemma \ref{dia} (1). We come to a contradiction since $0.55+0.05<0.85.$
\endproof

\begin{lemma}\label{quasi} Assume that $W$ is a periodic word with minimal period $A$ and $W=V$ in $H$ for some word $V$.

(1) Then $|V|\ge 0.25|W|.$ In particular, $A$ has infinite order
in $H.$

(2) A power $A^n$ cannot be conjugate to a word of length
$\le 0.2n|A|$.
\end{lemma}

\proof (1)  We consider a reduced diagram $\Delta$ with contour $pq$, where
$Lab(p)\equiv W$ and $Lab(q^{-1})\equiv V$. By induction on $|W|$ one may assume that the
path $p$ has no common edges with $q^{-1}$ and with $p^{-1}$. It follows that every (non-oriented) edge of $p$
is shared with a cell of $\Delta$, and by lemmas \ref{1touch}(2) and \ref{0.55}, the sum $\Sigma$ of the reduced perimeters of all cells in $\Delta$ is at least $\frac{20}{11}|p|$. On the other hand,
by Lemma \ref{Grin} (2), the number of (non-oriented) inner edges in $\Delta$ is at most $0.15\Sigma$. Hence the path $q$
has to have at least $(1- \frac{11}{20}-2\times 0.15)\Sigma=0.15\Sigma>0.25|p|$ edges, as required.

(2) Using the same argument, we have $|V|>0.25 |A^m|$ if $V=A^m$ in $H$ and $m\ge 1$. Now if $U=XA^nX^{-1}$ in $H$ and $|U|\le 0.2n|A|,$ then  $|U^s|\le 0.2sn|A|$ for every $s>0$, whence $A^{sn}$ is equal to the word $X^{-1}U^sX$ having length
 at most $0.2sn|A|+2|X|$. But this sum is less then $0.25sn|A|$ if $s$ is large enough, a contradiction.

\endproof

We will denote by $q_-$ (by $q_+$) the initial (resp., terminal) vertex of a path $q$.

 \begin{lemma} \label{polosa} Let $A$ be a minimal word and $\Delta$ be a reduced
diagram with contour $p_1q_1p_2q_2$, where $|p_1|+|p_2|\le 0.01(|q_1|+|q_2|)$
and $Lab(q_1)$, $Lab(q_2)$ are $A^{\pm 1}$-periodic words. Then

(1) there is a path $x$ of length $<27|A|$ connecting in $\Delta$ a vertex $o$ of $q_1$,
$o\ne (q_1)_-$, with a vertex of $q_2$;

(2) assume that $q_1=q'q''$, where $|q'|\ge 100(|p_1|+27|A|)$ and $|q''|\ge 100(|p_2|+27|A|)$. Then the vertex $q'_+$
can be connected with a vertex of $q_2$ by a path of length $<2800|A|$.
\end{lemma}

\proof (1) It is nothing to prove if $q_1$ and $q_2^{-1}$ share an edge. So we may assume that
they have no edges in common.

By Lemmas \ref{1touch} (2) and \ref{0.55}, there is no cell $\Pi$ in $\Delta$ such that $\ge 0.55|\partial'\Pi|$
of its boundary edges belong to some $q_i$, $i=1,2$.
Therefore by Lemma \ref{obe}, there exists a cell $\Pi$ having $>0.55|\partial'\Pi|$ edges on $q_1$ and $q_2$
but not on the one of these paths.
Hence the reduced contour of $\Pi$ factorizes as
$x_1y_1x_2y_2$ where $y_1$ and $y_2$ start and end on $q_1$ and $q_2$, resp., and $|x_1|, |x_2|<0.45|\partial'\Pi|$.
One of the two
paths $y_1$, $y_2$, call it $y$, has length at least $0.55|\partial'\Pi|/2>0.27|\partial'\Pi|$. Also $|y|<16 |A|$ by Lemma \ref{1.1} (2).
 Therefore $|x_1|, |x_2|<\frac{45}{27}|y|<27|A|$.
Finally, $x$ is either $x_1$ or $x_2$.

(2) Let $o$ be the first vertex on $q_1$ which can be connected with a vertex of the path $q_2$
by a path $x$ of length $<27|A|$. Such a vertex exists by the statement (1). Moreover, if $q_1=yz$,
where $y_+=o$, then $|y|<100(|p_1|+27|A|)$ since otherwise one can apply statement (1) to
the subdiagram with the boundary $p_1yxy'$, where $y'$ is the subpath of $q_2$, and replace
$o$ by another vertex on $y$.

Similarly, one can find the last vertex $o'$ on $q_1$ connected with $q_2$ by a path $x'$ of length
$<27|A|$, and so the vertex $q'_+$ lies between $o$ and $o'$. One can draw $x$ and $x'$ so that they
do not cross each other. Consider now the subdiagram
with contour $x^{-1}\bar q_1x'\bar q_2$, where $\bar q_i$ is a subpath of $q_i$ ($i=1,2$),
and repeat the same trick cutting up this diagram and obtaining subdiagrams with
contours $x_i^{-1}s_ix_{i+1}t_i$, where $s_1\dots s_k$ is the subpath of $q_1$ starting at $o$
and ending at $o'$, $|x_i|<27|A|$ and $|s_i|\le 5400|A|$. Since the vertex $q'_+$ belongs
to some $s_i$, it can be connected with a vertex of $t_i$ (and of $q_2$) by a path of
length $<2700|A|+27|A|<2800|A|.$
\endproof

 \begin{lemma} \label{pnqn} Let $A$ be a  minimal  word and $\Delta$ be a reduced
diagram with  contour $p_1q_1\dots p_nq_n$, where the factors are reduced paths, the words $Lab(q_i)$ are $A$- or $A^{-1}$-periodic
($i=1,\dots,n$) of length $\ge 100|A|$ and $|p_i|\le 0.01|q_j|$ for all $i,j$. Then for some $i$, $\Delta$ has
a subdiagram $\Gamma$ with contour $q'p_iq''x$, where $q'$ and $q''$ are the subpaths of $q_{i-1}$ and $q_i$, resp. (indices are
taken modulo $n$), $|x|<15|A|$ and either  $|q'|>|q_{i-1}|/3$ or  $|q''|>|q_i|/3.$
\end{lemma}

\proof We will construct a nested series of subdiagrams $\Delta\supset\Delta_1\supset\dots$ using
 transformations of two types. The transformations (a) and (b) of the first type
 %$\Delta\to \Delta_1$
will affect only one of the 'long' subpaths $q_i$-s. The transformations (a) and (b) of the second type
will affect two (cyclically) neighbors $q_{i-1}$ and $q_i$ for some $i$.

{\bf Type 1.} (a) Assume that we have $p_i=p^1_ie$ and $q_i=e^{-1}q^1_i$ (or $p_i=ep_i^1$ and $q_{i-1}=q^1_{i-1}e^{-1}$), where $e$ is an edge. Then we can replace the  subpath $p_iq_i$ in $\partial\Delta$ by $p^1_iq^1_i$. The diagram $\Delta_1$
has contour $p^1_1q^1_1\dots p^1_np^1_n$, where $p^1_j=p_j$ and $q^1_j=q_j$ for $j\ne i$.

(b) Assume that the path $p_iq_i$ (or $q_{i-1}p_i$) is reduced and it contains a Grindlinger arc $p$ of a cell $\Pi$,
i.e., $px$ is the reduced contour of $\Pi$ with $|x|<0.15|\partial'\Pi|$.

If $p$ is a subpath of $p_i$, we just replace it by $x^{-1}$ in $p_i$ and obtain the shorter path $p_i^1$
in the diagram $\Delta_1=\Delta\backslash \Pi$. Similar transformation works in any case if $p$ and $p_i$ have
a common subpath of length $>0.5|\partial'\Pi|$.

Then we assume that $p=yz$, where $p_i=uy$, $q_i=zq^1_i$ (Fig. 3) and $|z|<0.55|\partial'\Pi|$ by Lemma \ref{0.55}.
It follows that $|y|>(0.85-0.55)|\partial'\Pi|=0.3|\partial'\Pi|$ and $|y|-|x|>0.15|\partial'\Pi|$.
Now we define $p^1_i=ux^{-1}$ and removing $\Pi$, obtain a diagram $\Delta_1$, where

\begin{equation}\label{otrez}
|q_{i-1}|-|q^1_{i-1}|< \frac{0.55}{0.15}(|p_i|-|p_i^1|)< 4(|p_i|-|p^1_i|)
\end{equation}

\begin{figure}[h!]
\begin{center}
% This is a LaTeX picture output by TeXCAD.
% File name: [emb2.pic].
% Version of TeXCAD: 4.3
% Reference / build: 30-Jun-2012 (rev. 105)
% For new versions, check: http://texcad.sf.net/
% Options on the following lines.
%\grade{\on}
%\emlines{\off}
%\epic{\off}
%\beziermacro{\on}
%\reduce{\on}
%\snapping{\off}
%\pvinsert{% Your \input, \def, etc. here}
%\quality{8.000}
%\graddiff{0.005}
%\snapasp{1}
%\zoom{4.0000}
\unitlength 1mm % = 2.845pt
\linethickness{0.4pt}
\ifx\plotpoint\undefined\newsavebox{\plotpoint}\fi % GNUPLOT compatibility
\begin{picture}(171.25,65.75)(12,0)
%\emline(10,52.25)(41,20.5)
\multiput(10,52.25)(.03373231774,-.0345484222){919}{\line(0,-1){.0345484222}}
%\end
%\emline(37.75,52)(29,32.5)
\multiput(37.75,52)(-.0336538462,-.075){260}{\line(0,-1){.075}}
%\end
\put(27.5,44.5){\makebox(0,0)[cc]{$\Pi$}}
\put(22,31){\makebox(0,0)[cc]{$p_i$}}
\put(36.75,27.5){\makebox(0,0)[cc]{$u$}}
\put(37.25,44){\makebox(0,0)[cc]{$x$}}
\put(18.75,46){\makebox(0,0)[cc]{$y$}}
\put(22.25,50.25){\makebox(0,0)[cc]{$z$}}
\put(43.25,57.5){\makebox(0,0)[cc]{$q_i$}}
\put(62,49){\makebox(0,0)[cc]{$q_i^1$}}
\put(58.5,36.25){\makebox(0,0)[cc]{$\Delta_1$}}
%\emline(30.5,53.25)(33.25,52.25)
\multiput(30.5,53.25)(.0916667,-.0333333){30}{\line(1,0){.0916667}}
%\end
%\emline(30.75,51.25)(33,52.25)
\multiput(30.75,51.25)(.075,.0333333){30}{\line(1,0){.075}}
%\end
%\emline(31.75,41.5)(31.25,38.75)
\multiput(31.75,41.5)(-.0333333,-.1833333){15}{\line(0,-1){.1833333}}
%\end
%\emline(31.25,38.75)(33.75,40.5)
\multiput(31.25,38.75)(.04807692,.03365385){52}{\line(1,0){.04807692}}
%\end
%\emline(19.5,42.5)(20.75,40)
\multiput(19.5,42.5)(.03289474,-.06578947){38}{\line(0,-1){.06578947}}
%\end
%\emline(19.5,42.5)(21.75,41.25)
\multiput(19.5,42.5)(.05921053,-.03289474){38}{\line(1,0){.05921053}}
%\end
%\emline(30.75,31)(32.25,28.5)
\multiput(30.75,31)(.03333333,-.05555556){45}{\line(0,-1){.05555556}}
%\end
%\emline(31.5,30.5)(33,29.5)
\multiput(31.5,30.5)(.05,-.0333333){30}{\line(1,0){.05}}
%\end
%\emline(69.75,53.75)(72.25,52.75)
\multiput(69.75,53.75)(.0833333,-.0333333){30}{\line(1,0){.0833333}}
%\end
%\emline(70,51)(72.25,52.75)
\multiput(70,51)(.04326923,.03365385){52}{\line(1,0){.04326923}}
%\end
\put(35.25,52.25){\line(0,-1){4.75}}
\put(33.25,52.25){\line(0,-1){10.25}}
\put(31,52.5){\line(0,-1){13.5}}
\put(28.75,52.25){\line(0,-1){3.75}}
\put(28.75,42.75){\line(0,-1){8.75}}
\put(27,52.25){\line(0,-1){2.5}}
\put(27,42.75){\line(0,-1){6.75}}
\put(25,51.75){\line(0,-1){3}}
\put(25.25,42.75){\line(0,-1){5.75}}
\put(23.25,51.75){\line(0,-1){12}}
\put(21.25,49.75){\line(-1,0){.25}}
\put(21.5,49.5){\line(0,-1){8}}
\put(19.75,52){\line(0,-1){9.75}}
%\emline(17.75,52.25)(18,49.25)
\multiput(17.75,52.25)(.03125,-.375){8}{\line(0,-1){.375}}
%\end
\put(16,52.25){\line(0,-1){5.25}}
%\emline(14,52.25)(14.25,52.5)
\multiput(14,52.25)(.03125,.03125){8}{\line(0,1){.03125}}
%\end
\put(14.25,51.75){\line(0,-1){3.25}}
\put(12.25,52){\line(0,-1){1.75}}
%\emline(97.25,51.5)(97.5,51.75)
\multiput(97.25,51.5)(.03125,.03125){8}{\line(0,1){.03125}}
%\end
\thicklines
\put(97.5,51.25){\line(1,0){69.75}}
\put(97.75,51.5){\line(0,-1){31}}
\put(98,20.5){\line(1,0){69.5}}
\put(127,51){\line(0,-1){29.75}}
\put(91.5,36){\makebox(0,0)[cc]{$p_i$}}
\put(132,15.75){\makebox(0,0)[cc]{$q_{i-1}$}}
\put(131,55){\makebox(0,0)[cc]{$q_i$}}
\put(153,24.75){\makebox(0,0)[cc]{$q_{i-1}^1$}}
\put(152.25,47){\makebox(0,0)[cc]{$q_i^1$}}
\put(161.25,36.5){\makebox(0,0)[cc]{$\Delta_1$}}
\put(112,47.75){\makebox(0,0)[cc]{$z$}}
\put(112,23.25){\makebox(0,0)[cc]{$y$}}
\put(110.5,35.5){\makebox(0,0)[cc]{$\Pi$}}
\put(129.75,38.75){\makebox(0,0)[cc]{$x$}}
%\emline(98,45)(96.75,42.75)
\multiput(98,45)(-.03289474,-.05921053){38}{\line(0,-1){.05921053}}
%\end
%\emline(98.25,44.75)(98.75,43.5)
\multiput(98.25,44.75)(.0333333,-.0833333){15}{\line(0,-1){.0833333}}
%\end
%\emline(103.75,20.5)(106,21.5)
\multiput(103.75,20.5)(.075,.0333333){30}{\line(1,0){.075}}
%\end
%\emline(104,20.75)(105.5,20)
\multiput(104,20.75)(.0652174,-.0326087){23}{\line(1,0){.0652174}}
%\end
\put(104.5,20.5){\line(1,0){.25}}
%\emline(104.25,20.5)(106,19.75)
\multiput(104.25,20.5)(.076087,-.0326087){23}{\line(1,0){.076087}}
%\end
%\emline(116.75,52.25)(119,51.25)
\multiput(116.75,52.25)(.075,-.0333333){30}{\line(1,0){.075}}
%\end
%\emline(117.25,50.5)(118,51)
\multiput(117.25,50.5)(.05,.0333333){15}{\line(1,0){.05}}
%\end
%\emline(116.75,50.5)(118.5,51.75)
\multiput(116.75,50.5)(.04605263,.03289474){38}{\line(1,0){.04605263}}
%\end
%\emline(126,31.25)(126.5,29.75)
\multiput(126,31.25)(.0333333,-.1){15}{\line(0,-1){.1}}
%\end
%\emline(128,31.25)(127.25,30)
\multiput(128,31.25)(-.0326087,-.0543478){23}{\line(0,-1){.0543478}}
%\end
%\emline(159.75,52.5)(162.5,51.5)
\multiput(159.75,52.5)(.0916667,-.0333333){30}{\line(1,0){.0916667}}
%\end
%\emline(160,50.25)(162,51.5)
\multiput(160,50.25)(.05263158,.03289474){38}{\line(1,0){.05263158}}
%\end
%\emline(159.75,21)(162.5,21.75)
\multiput(159.75,21)(.1195652,.0326087){23}{\line(1,0){.1195652}}
%\end
%\emline(160,20.25)(162.25,19.5)
\multiput(160,20.25)(.0978261,-.0326087){23}{\line(1,0){.0978261}}
%\end
\put(371.25,99.75){\line(-1,0){.5}}
\thinlines
%\emline(98,47)(103.25,51.5)
\multiput(98,47)(.039179104,.03358209){134}{\line(1,0){.039179104}}
%\end
%\emline(98.25,43.5)(107.25,51.5)
\multiput(98.25,43.5)(.037815126,.033613445){238}{\line(1,0){.037815126}}
%\end
%\emline(98,40.5)(110.25,51)
\multiput(98,40.5)(.0392628205,.0336538462){312}{\line(1,0){.0392628205}}
%\end
%\emline(98,37.75)(110.5,48.25)
\multiput(98,37.75)(.0400641026,.0336538462){312}{\line(1,0){.0400641026}}
%\end
%\emline(98,34.5)(117,50.5)
\multiput(98,34.5)(.04,.0336842105){475}{\line(1,0){.04}}
%\end
%\emline(98,31.75)(120.75,50.75)
\multiput(98,31.75)(.0403368794,.0336879433){564}{\line(1,0){.0403368794}}
%\end
\put(98.25,28.75){\line(6,5){9}}
%\emline(111.75,39.5)(124,50.75)
\multiput(111.75,39.5)(.0366766467,.0336826347){334}{\line(1,0){.0366766467}}
%\end
%\emline(98,25.5)(107.25,33.5)
\multiput(98,25.5)(.038865546,.033613445){238}{\line(1,0){.038865546}}
%\end
%\emline(113.75,38.75)(127.25,50.75)
\multiput(113.75,38.75)(.0379213483,.0337078652){356}{\line(1,0){.0379213483}}
%\end
%\emline(97.75,22.5)(126.75,47)
\multiput(97.75,22.5)(.03988995873,.03370013755){727}{\line(1,0){.03988995873}}
%\end
%\emline(100,21.25)(126.25,43.75)
\multiput(100,21.25)(.03935532234,.03373313343){667}{\line(1,0){.03935532234}}
%\end
%\emline(103,21)(126.75,41.25)
\multiput(103,21)(.0395174709,.0336938436){601}{\line(1,0){.0395174709}}
%\end
%\emline(106.75,21.25)(110.25,25)
\multiput(106.75,21.25)(.033653846,.036057692){104}{\line(0,1){.036057692}}
%\end
%\emline(110.25,25)(126,38)
\multiput(110.25,25)(.0408031088,.0336787565){386}{\line(1,0){.0408031088}}
%\end
%\emline(110,21.5)(126.5,36.25)
\multiput(110,21.5)(.0376712329,.0336757991){438}{\line(1,0){.0376712329}}
%\end
%\emline(112.75,21.25)(126.25,33.5)
\multiput(112.75,21.25)(.0370879121,.0336538462){364}{\line(1,0){.0370879121}}
%\end
%\emline(115.75,21)(126.75,31.5)
\multiput(115.75,21)(.0352564103,.0336538462){312}{\line(1,0){.0352564103}}
%\end
%\emline(119,21.25)(127,28.75)
\multiput(119,21.25)(.035874439,.033632287){223}{\line(1,0){.035874439}}
%\end
\put(121.5,20.75){\line(1,1){5.5}}
%\emline(124.5,21.25)(126.5,23.75)
\multiput(124.5,21.25)(.03333333,.04166667){60}{\line(0,1){.04166667}}
%\end
\put(48.5,10.25){\makebox(0,0)[cc]{Figure 3}}
\put(142,11.25){\makebox(0,0)[cc]{Figure 4}}
\thicklines
%\emline(10.25,52.5)(75.5,52.75)
\multiput(10.25,52.5)(8.15625,.03125){8}{\line(1,0){8.15625}}
%\end
\end{picture}

\end{center}
\end{figure}

{\bf Type 2.} (a) Assume that $|p_i|=0$ for some $i$ and $q_{i-1}=q_{i-1}^1e$, $q_i=e^{-1}q^1_i$ for an edge $e$.
Then we replace the subpath $q_{i-1}p_iq_i$ by $q^1_{i-1}p^1_iq^1_i$, where $|p^1_i|=0$.

(b) Let for some $i$, we have a Grindlinger cell $\Pi$ with the Grindlinger arc
 $p=yp_iz$, where $|p_i|\le  |\partial'\Pi|/2$ and $q_{i-1}=q^1_{i-1}y$ and $q_i=zq^1_i$ (Fig. 4).
 Then we have $|y|+|z|>(0.85-0.5)|\partial'\Pi|=
0.35|\partial'\Pi|$. But $|y|, |z|\le 16|A|$ by Lemma \ref{1.1} (2), and so $32|A|>0.35|\partial\Pi|$,
whence $|\partial'\Pi|<100|A|$ and for the compliment $x$ of $p$ in $\partial'\Pi$, we have
$|x|<0.15\times 100|A|=15|A|.$

In this case, we define the contour of $\Delta_1$ ($=\Delta\backslash\Pi$) as follows.
we replace $p_i$ by $p^1_i=x^{-1}$ and replace $q_{i-1}$ and $q_i$ by $q^1_{i-1}$ and $q^1_i$, resp.

\begin{rk} \label{pq} One may assume that the paths $q^1_{i-1}x^{-1}$ and $x^{-1}q^1_i$ are reduced after
the transformation of type 2(b) since the Grindlinger arc $p$ can be made the longest possible.
Therefore if  $\Delta_1$ results from $\Delta$ after a transformation of the second type
at the 'corner' $i$, then no transformation of the first type is applicable to the same $i$-th corner
of $\Delta_1$. Indeed, the path $x$ was a part of the boundary of $\Pi$. Hence the boundary of a cell $\pi$ from $\Delta_1$ cannot share with $x$ a path of length $\ge 0.05|\partial\pi|$ by Lemma \ref{dia} (1). Since $\partial'\pi$ cannot share
a path of length $\ge 0.55|\partial'\pi|$ with $q_{i-1}^1$ or with $q_i^1$ by Lemma \ref{0.55}, it cannot be a Grindlinger cell
of type 1 in $\Delta_1$ because $0.05+0.55<0.85$.
\end{rk}

The transformations of the first type $\Delta\to\Delta_1\to\Delta_2\dots$ can shorten the path $q_i$
by at most $4(|p_i|+|p_{i+1}|)$, as it follows from (\ref{otrez}). Therefore after a maximal chain of the transformations
of the first type applied to all the corners $i=1,\dots,n$ in arbitrary order, we will have a subdiagram $\bar\Delta$ with contour
$\bar p_1\bar q_1\dots \bar p_n\bar q_n$, such that $|\bar q_i|>|q_i|-0.08|q_i|> 0.9|q_i|$, $|\bar p_i|\le|p_i|$
for every $i=1,\dots,n.$

By Remark \ref{pq}, only transformations of the second type can appear in any decreasing chain
of subdiagrams $\bar\Delta\to\bar\Delta_1\to\dots$. If $\bar\Delta_j$ has contour
$\bar p^j_1\bar q^j_1\dots \bar p^j_n\bar q^j_n$ and for each $i$ we have $|\bar q^j_i|\ge |q_i|/3>33|A|$, then one more
transformation of the second type decreases the perimeter of $\Delta_j$, because
no Grindlinger arc of a cell  can entirely cover some $\bar q^j_i$ by Lemma \ref{1.1} (2).
By the same lemma, any
such transformation affects only an initial or a terminal part of length $<16|A|$ for
some $\bar q^j_{i-1}$ and $\bar q^j_i.$ Thus, after some transformation of the second type we
can obtain $|\bar q^j_i|<|q_i|/3$, and so $q_i=q''\bar q^j_iq'''$, where $|q''|>|q_i|/3$ or $|q'''|>|q_i|/3$.
Let us consider only the first event. Then the vertex $q''_+$ is connected with $q_{i-1}$ by a path $x$
of length $<15|A|$ since this $x$ appears after a transformation of the second type. The
desired subdiagram $\Gamma$ is obtained.
\endproof

Below  a reduced disc diagram $\Delta$ is  called a {\it box} if its contour is a product
%of four reduced subpaths
$x_1y_1x_2y_2$,
where the subpaths $x_1, y_1, y_2$ are reduced, $x_2$ is reduced or a $G$-fragment,
$Lab(x_1)$ is an $A$-periodic word with a minimal period $A$, and $y_1^{-1}$, $y_2$ are {\it perpendicular}
to $x_1$, i.e., the distance in $\Delta$ between $(y_1)_+$ and $x_1$ (between $(y_2)_-$ and $x_1$) is
equal to $|y_1|$ (to $|y_2|$, resp.).
%, and (3) $|x_2|\le c (|y_1|+|y_2|)$ for $c=$.

A subdiagram $\Delta'$ of a box $\Delta$ with a contour $x'_1y'_1x'_2y'_2$ is a {\it subbox} if $x'_1=x_1$, $y'_1$ (resp., $y'_2$)
is the beginning of the path $y_1$ (the end of $y_2$, and $x'_2$ is a reduced path or a $G$-fragment. So a subbox is a box itself.

\begin{lemma} \label{frag} Let $\Delta$ be a box with $|x_1|\ge 16|A|$. Then $|x_2|>1$,
and $x_2$ cannot be a $G$-fragment.
%if $\Delta$ has no cells compatible with $x_2$.
\end{lemma}
\proof Assume on the contrary that $|x_2|\le 1$ or
%\Delta$ has a cell compatible with the $G$-fragment
 $x_2$ is a $G$-gragment. Then $\Delta$ has a subbox $\Delta'$, which is minimal
with respect to the property that either $|x'_2|\le 1$ or $x'_2$ is a $G$-fragment.

Note that  the subpaths $x'_1y'_1$ and $y'_2x'_1$
are reduced since $(y'_1)^{-1}$ and $y'_2$ are perpendicular to $x'_1$, while $y'_1x'_2$ and $x'_2y'_2$ (or
$y'_1y'_2$ if $|x'_2|=0$) are reduced due the minimality of the subbox $\Delta'$.
Also one cannot have $|y'_1|=|y'_2|=0$ since $Lab(x_1)=Lab(x'_1)\ne Lab(x'_2)$ in the free group by Lemma \ref{1.1} (2).
Therefore the boundary $\partial\Delta'$ is not a completely cancellable path, whence $\Delta'$ must have $\cal R$-cells, and so it has a Grindlinger cell $\pi$ by Lemma \ref{Grin} (1).

The Grindlinger arc $p$ of $\partial'\pi$ cannot have common edges with both paths $y'_1$ and $y'_2$
since in this case one obtains a smaller subbox $\Delta''=\Delta'\backslash\pi$,
where $x''_2$ is a $G$-fragment (Fig. 5).
 Also $p$ does not include a part of $x'_2$ of length $\ge 0.05|\partial'\pi|$
since otherwise one can obtain a smaller counter-example using  Lemma \ref{dia} (2).
 The path $p$ does not include a subpath of $y'_1$ (or of $y'_2$) with length $>0.5|\partial'\pi|$ since the perpendicular $(y'_1)^{-1}$ is geodesic in $\Delta$. Finally,
the arc $p$ cannot cover the whole path $x'_1$ by Lemma \ref{1.1} (2), and the subpath of $p$ shared with
$(x'_1)^{\pm 1}$ has length $<0.55|\partial'\pi|$ by Lemma \ref{0.55}.

\begin{figure}[h!]
\begin{center}
% This is a LaTeX picture output by TeXCAD.
% File name: [emb3.pic].
% Version of TeXCAD: 4.3
% Reference / build: 30-Jun-2012 (rev. 105)
% For new versions, check: http://texcad.sf.net/
% Options on the following lines.
%\grade{\on}
%\emlines{\off}
%\epic{\off}
%\beziermacro{\on}
%\reduce{\on}
%\snapping{\off}
%\pvinsert{% Your \input, \def, etc. here}
%\quality{8.000}
%\graddiff{0.005}
%\snapasp{1}
%\zoom{4.0000}
\unitlength 1mm % = 2.845pt
\linethickness{0.4pt}
\ifx\plotpoint\undefined\newsavebox{\plotpoint}\fi % GNUPLOT compatibility
\begin{picture}(200.5,67.25)(15,0)
\thicklines
%\emline(100,24.25)(139.75,56)
\multiput(100,24.25)(.04219745223,.03370488323){942}{\line(1,0){.04219745223}}
%\end
\put(139.75,56){\line(1,0){37.25}}
\put(177,56){\line(0,-1){33.25}}
%\emline(177,22.75)(99,23)
\multiput(177,22.75)(-9.75,.03125){8}{\line(-1,0){9.75}}
%\end
%\emline(101,24.5)(98.75,23)
\multiput(101,24.5)(-.05,-.03333333){45}{\line(-1,0){.05}}
%\end
%\emline(128.25,47)(138.5,23.5)
\multiput(128.25,47)(.0337171053,-.0773026316){304}{\line(0,-1){.0773026316}}
%\end
%\emline(142.75,61)(139.5,56.5)
\multiput(142.75,61)(-.033505155,-.046391753){97}{\line(0,-1){.046391753}}
%\end
\put(176.75,60){\line(0,-1){4.75}}
%\emline(175.75,41.25)(176.75,37)
\multiput(175.75,41.25)(.0333333,-.1416667){30}{\line(0,-1){.1416667}}
%\end
%\emline(178.25,41.25)(177,38)
\multiput(178.25,41.25)(-.03289474,-.08552632){38}{\line(0,-1){.08552632}}
%\end
%\emline(167.5,57.5)(171,56.25)
\multiput(167.5,57.5)(.09210526,-.03289474){38}{\line(1,0){.09210526}}
%\end
%\emline(167.25,54.75)(170.75,56)
\multiput(167.25,54.75)(.09210526,.03289474){38}{\line(1,0){.09210526}}
%\end
%\emline(116.25,23.25)(120.25,24.5)
\multiput(116.25,23.25)(.10526316,.03289474){38}{\line(1,0){.10526316}}
%\end
%\emline(116.5,23.25)(120,21.5)
\multiput(116.5,23.25)(.06730769,-.03365385){52}{\line(1,0){.06730769}}
%\end
%\emline(117.75,39.5)(121.5,41.5)
\multiput(117.75,39.5)(.0625,.03333333){60}{\line(1,0){.0625}}
%\end
%\emline(121.5,41.5)(119.25,38.5)
\multiput(121.5,41.5)(-.03358209,-.04477612){67}{\line(0,-1){.04477612}}
%\end
%\emline(133.25,32.25)(135.5,30.25)
\multiput(133.25,32.25)(.0375,-.03333333){60}{\line(1,0){.0375}}
%\end
\put(135.25,33.75){\line(0,-1){2.75}}
\put(126,36.25){\makebox(0,0)[cc]{$\pi$}}
\put(150,44.75){\makebox(0,0)[cc]{$\Delta'$}}
\put(180.5,46){\makebox(0,0)[cc]{$y'_2$}}
\put(154.25,61){\makebox(0,0)[cc]{$x'_2$}}
\put(143.25,18.25){\makebox(0,0)[cc]{$x'_1=x_1$}}
\put(124.25,26.75){\makebox(0,0)[cc]{$u$}}
\put(117,33.25){\makebox(0,0)[cc]{$v$}}
\put(120.25,47.75){\makebox(0,0)[cc]{$y'_1$}}
\put(134.5,38.5){\makebox(0,0)[cc]{$q$}}
%\emline(8.5,24.5)(29.75,58.25)
\multiput(8.5,24.5)(.0337301587,.0535714286){630}{\line(0,1){.0535714286}}
%\end
\put(29.75,58.25){\line(1,0){27.75}}
%\emline(57.5,58.25)(71.25,24.25)
\multiput(57.5,58.25)(.0337009804,-.0833333333){408}{\line(0,-1){.0833333333}}
%\end
\put(71.25,24.25){\line(-1,0){63.25}}
%\emline(27,54)(32.5,47.5)
\multiput(27,54)(.033536585,-.039634146){164}{\line(0,-1){.039634146}}
%\end
\put(32.5,47.5){\line(1,0){19.75}}
%\emline(52.25,47.5)(59.25,55)
\multiput(52.25,47.5)(.033653846,.036057692){208}{\line(0,1){.036057692}}
%\end
\put(13.75,44){\makebox(0,0)[cc]{$y'_1$}}
\put(66.75,44){\makebox(0,0)[cc]{$y'_2$}}
\put(31,19.5){\makebox(0,0)[cc]{$x''_1=x'_1=x_1$}}
\put(43.5,31.5){\makebox(0,0)[cc]{$\Delta''$}}
\put(42.75,62.25){\makebox(0,0)[cc]{$x'_2$}}
\put(43.75,44){\makebox(0,0)[cc]{$x''_2$}}
\put(40.75,52.75){\makebox(0,0)[cc]{$\pi$}}
\thinlines
\put(32.75,57.75){\line(-2,-3){3.5}}
%\emline(36,58)(30.5,49)
\multiput(36,58)(-.033536585,-.054878049){164}{\line(0,-1){.054878049}}
%\end
%\emline(38.75,58)(33,48)
\multiput(38.75,58)(-.033625731,-.058479532){171}{\line(0,-1){.058479532}}
%\end
%\emline(42,58)(40.5,56)
\multiput(42,58)(-.03333333,-.04444444){45}{\line(0,-1){.04444444}}
%\end
%\emline(38.25,52.5)(35.75,48)
\multiput(38.25,52.5)(-.03333333,-.06){75}{\line(0,-1){.06}}
%\end
%\emline(44.5,58.25)(42.5,56.25)
\multiput(44.5,58.25)(-.03333333,-.03333333){60}{\line(0,-1){.03333333}}
%\end
%\emline(40.25,51.25)(38.25,48.5)
\multiput(40.25,51.25)(-.03333333,-.04583333){60}{\line(0,-1){.04583333}}
%\end
%\emline(47,58.25)(40.5,48)
\multiput(47,58.25)(-.033678756,-.053108808){193}{\line(0,-1){.053108808}}
%\end
%\emline(49.75,58.5)(43.75,48.25)
\multiput(49.75,58.5)(-.033707865,-.05758427){178}{\line(0,-1){.05758427}}
%\end
%\emline(52.25,58.5)(46.25,48.25)
\multiput(52.25,58.5)(-.033707865,-.05758427){178}{\line(0,-1){.05758427}}
%\end
%\emline(54.75,58.75)(49,48.25)
\multiput(54.75,58.75)(-.033625731,-.061403509){171}{\line(0,-1){.061403509}}
%\end
%\emline(57,58.25)(51.75,48.75)
\multiput(57,58.25)(-.033653846,-.060897436){156}{\line(0,-1){.060897436}}
%\end
%\emline(58.25,55.5)(56,52)
\multiput(58.25,55.5)(-.03358209,-.05223881){67}{\line(0,-1){.05223881}}
%\end
\thicklines
%\emline(13,33.25)(16.75,37.25)
\multiput(13,33.25)(.033482143,.035714286){112}{\line(0,1){.035714286}}
%\end
%\emline(16.75,37.25)(14.75,32.75)
\multiput(16.75,37.25)(-.03333333,-.075){60}{\line(0,-1){.075}}
%\end
%\emline(58.75,26)(54.5,24.25)
\multiput(58.75,26)(-.08173077,-.03365385){52}{\line(-1,0){.08173077}}
%\end
%\emline(54.5,24.25)(58.25,23.5)
\multiput(54.5,24.25)(.1630435,-.0326087){23}{\line(1,0){.1630435}}
%\end
\put(68,36.25){\line(1,0){.25}}
%\emline(68.25,36.25)(68,33.25)
\multiput(68.25,36.25)(-.03125,-.375){8}{\line(0,-1){.375}}
%\end
%\emline(65.5,35.25)(67.25,33.25)
\multiput(65.5,35.25)(.03365385,-.03846154){52}{\line(0,-1){.03846154}}
%\end
%\emline(47,59.5)(50,58.5)
\multiput(47,59.5)(.1,-.0333333){30}{\line(1,0){.1}}
%\end
%\emline(47.5,57.25)(49.5,59)
\multiput(47.5,57.25)(.03846154,.03365385){52}{\line(1,0){.03846154}}
%\end
\put(49.5,59){\line(-1,0){.25}}
%\emline(47,57.25)(49.75,58)
\multiput(47,57.25)(.1195652,.0326087){23}{\line(1,0){.1195652}}
%\end
\put(33.5,10.75){\makebox(0,0)[cc]{Figure 5}}
\put(125.75,11.5){\makebox(0,0)[cc]{Figure 6}}
\thinlines
%\emline(109.25,32.25)(103.5,23.75)
\multiput(109.25,32.25)(-.033625731,-.049707602){171}{\line(0,-1){.049707602}}
%\end
%\emline(113.75,35.25)(106.5,24)
\multiput(113.75,35.25)(-.03372093,-.052325581){215}{\line(0,-1){.052325581}}
%\end
%\emline(115,32.75)(109.75,23.75)
\multiput(115,32.75)(-.033653846,-.057692308){156}{\line(0,-1){.057692308}}
%\end
%\emline(122.5,41.75)(112.25,23)
\multiput(122.5,41.75)(-.0337171053,-.0616776316){304}{\line(0,-1){.0616776316}}
%\end
%\emline(126,45.5)(115.5,23.5)
\multiput(126,45.5)(-.0336538462,-.0705128205){312}{\line(0,-1){.0705128205}}
%\end
%\emline(128.75,46.75)(126.25,40.75)
\multiput(128.75,46.75)(-.03333333,-.08){75}{\line(0,-1){.08}}
%\end
%\emline(123.25,35.5)(118,23.5)
\multiput(123.25,35.5)(-.033653846,-.076923077){156}{\line(0,-1){.076923077}}
%\end
%\emline(129.75,43.75)(128.25,40)
\multiput(129.75,43.75)(-.03333333,-.08333333){45}{\line(0,-1){.08333333}}
%\end
%\emline(125,34.5)(120.25,24.25)
\multiput(125,34.5)(-.033687943,-.072695035){141}{\line(0,-1){.072695035}}
%\end
\put(130.75,40.5){\line(-1,-2){5}}
%\emline(123.25,25.5)(122.25,23.75)
\multiput(123.25,25.5)(-.0333333,-.0583333){30}{\line(0,-1){.0583333}}
%\end
\put(131.75,38.25){\line(-1,-2){7}}
%\emline(133,35.75)(126.75,24.5)
\multiput(133,35.75)(-.033602151,-.060483871){186}{\line(0,-1){.060483871}}
%\end
%\emline(134.5,34)(128.75,24.5)
\multiput(134.5,34)(-.033625731,-.055555556){171}{\line(0,-1){.055555556}}
%\end
%\emline(135.25,31)(130.75,23.75)
\multiput(135.25,31)(-.03358209,-.054104478){134}{\line(0,-1){.054104478}}
%\end
%\emline(136.25,28.5)(132.5,23.25)
\multiput(136.25,28.5)(-.033482143,-.046875){112}{\line(0,-1){.046875}}
%\end
%\emline(136.75,26.25)(135.25,23.25)
\multiput(136.75,26.25)(-.03333333,-.06666667){45}{\line(0,-1){.06666667}}
%\end
\end{picture}

\end{center}
\end{figure}

Thus, up to substitution $x'_1 \leftrightarrow x'_2$, the arc $p$ has a subpath $uv$ of length $>(0.85-0.05)|\partial'\pi|$,
where $u$ is a subpath of       $x'_1$, $v$ is a subpath of $y'_1$ and $|v|>(0.8-0.55)|\partial'\pi|=0.25|\partial'\pi|$ (Fig. 6). Hence we have that $\partial'\pi=pq$, where $|q|< (1-0.8)|\partial'\pi|=0.2|\partial'\pi|$, and $q$ connects
the vertex $v_+$ with $x'_1$. This contradicts the assumption that $y_1^{-1}$ is perpendicular to $x_1$
because $|q|<|v|.$ The lemma is proved.
\endproof

\begin{lemma} \label{16} Let $\Delta$ be a box with $|x_1|\ge 16|A|$. Then $|x_1|+|x_2|>0.01(|y_1|+|y_2|)$.
\end{lemma}
\proof
Assume that $|x_1|+|x_2|\le 0.01(|y_1|+|y_2|)$. Then we can apply Lemma \ref{obe} to $\Delta$.
Indeed, (1) the paths $y_1$ and $y_2^{-1}$ have disjoint sets of edges because otherwise we should
have a subbox with contour $x'_1y'_1x'_2y'_2$, where $|x'_2|=0$, a contradiction
with Lemma \ref{frag}, and (2) no reduced boundary $\partial'\Pi$ has $>0.5|\partial'\Pi|$ edges
on $y_i$ ($i=1,2$) since the perpendiculars are geodesic paths.

So there is a cell $\Pi$ whose boundary has common edges with both $y'_1$ an $y'_2$.
But an arc of $\partial'\Pi$ connecting some vertices of $y'_1$ and $y'_2$ is labeled by
a $G$-fragment, which by Lemma \ref{frag}, gives an impossible subbox, a contradiction.
\endproof

\begin{lemma} \label{sq}  Let $A$ be a minimal word and $\theta\in (0, 1)$. Then there is a
positive number $c=c(\theta)$ with the following property. Let $\Delta$ be a reduced diagram with
a contour $pq$ such that $Lab(p)\equiv A^{2n}$ and the distance in $\Delta$ from the middle vertex $o$
of $p$ to the path $q$ is at least $\theta n|A|$. Then $|q|\ge cn^2.$
\end{lemma}

\proof Let $v_0,v_1,\dots$ be all the vertices of $q$ counted from $q_-$ to $q_+.$ We consider
the set of paths $y^0,y^1,\dots$, where $y^i$ is a perpendicular to $p$ drawn in $\Delta$ from  $v_i$.
One may assume that the different $y^i$ and $y^j$ do not cross each other. Indeed, if they cross at a
vertex $u$ and $y^i=(y')^i(y'')^i$, $y^j=(y')^j(y'')^j$ with $(y')^i_+=(y')^j_+=u,$ then $|(y'')^i|=|(y'')^j|$
since both these subpaths are the shortest ones between $u$ and $p$, and so one can replace $y^j$
by $(y')^j(y'')^i$, and so on. Hence the ends $y^0_+, y^1_+,\dots$ are placed on $p$ from $p_+$ to $p_-$.

Let $\Delta_i$ be the box with contour $x_1^i(y^i)^{-1}x_2^iy^{i+1}$, where $x_1^i$ is a subpath of $p$
connecting $y^{i+1}_+$ and $y^i_+$, and $x_2^i$ is the path of length $1$ connecting $v_i$ and $v_{i+1}$,
$i=0,1, \dots$ (Fig. 7). We have $p=\dots x^1_1x^0_1$,  where for every $i,$ we obtain the inequality $|x^i_1|<16|A| $
from Lemma \ref{frag} since $|x^i_2|=1.$ So it is possible to find a pair of indices $(i,j)$ such that
the path $z=x_1^i\dots x_1^j$ contains the vertex $o$, $|z|>\theta n|A| -32|A|$, and the distances
between $o$ and the vertices $z_-$ and $z_+$ are less then $\theta n|A|/2$. (Here and below one may assume
that $n$ is large enough in comparison with $\theta^{-1}|A|$.)

\begin{figure}
\begin{center}

% This is a LaTeX picture output by TeXCAD.
% File name: [emv4.pic].
% Version of TeXCAD: 4.3
% Reference / build: 30-Jun-2012 (rev. 105)
% For new versions, check: http://texcad.sf.net/
% Options on the following lines.
%\grade{\on}
%\emlines{\off}
%\epic{\off}
%\beziermacro{\on}
%\reduce{\on}
%\snapping{\off}
%\pvinsert{% Your \input, \def, etc. here}
%\quality{8.000}
%\graddiff{0.005}
%\snapasp{1}
%\zoom{4.0000}
\unitlength 1mm % = 2.845pt
\linethickness{0.4pt}
\ifx\plotpoint\undefined\newsavebox{\plotpoint}\fi % GNUPLOT compatibility
\begin{picture}(162.25,44.5)(7,0)
\thicklines
\put(10.25,15.5){\line(1,0){61.75}}
\put(100.5,16.25){\line(1,0){61.75}}
%\emline(10.25,15)(20.25,31)
\multiput(10.25,15)(.0336700337,.0538720539){297}{\line(0,1){.0538720539}}
%\end
%\emline(100.5,15.75)(110.5,31.75)
\multiput(100.5,15.75)(.0336700337,.0538720539){297}{\line(0,1){.0538720539}}
%\end
%\emline(20.25,31)(37.25,40.25)
\multiput(20.25,31)(.0618181818,.0336363636){275}{\line(1,0){.0618181818}}
%\end
%\emline(110.5,31.75)(127.5,41)
\multiput(110.5,31.75)(.0618181818,.0336363636){275}{\line(1,0){.0618181818}}
%\end
\put(37.25,40.25){\line(1,0){24.75}}
\put(127.5,41){\line(1,0){24.75}}
%\emline(62.25,40.25)(69.5,29.5)
\multiput(62.25,40.25)(.03372093,-.05){215}{\line(0,-1){.05}}
%\end
%\emline(152.5,41)(159.75,30.25)
\multiput(152.5,41)(.03372093,-.05){215}{\line(0,-1){.05}}
%\end
%\emline(69.5,29.5)(71.75,15.75)
\multiput(69.5,29.5)(.03358209,-.20522388){67}{\line(0,-1){.20522388}}
%\end
%\emline(159.75,30.25)(162,16.5)
\multiput(159.75,30.25)(.03358209,-.20522388){67}{\line(0,-1){.20522388}}
%\end
%\emline(56.5,40.5)(54,15)
\multiput(56.5,40.5)(-.03333333,-.34){75}{\line(0,-1){.34}}
%\end
%\emline(62,40)(63.25,15)
\multiput(62,40)(.03289474,-.65789474){38}{\line(0,-1){.65789474}}
%\end
%\emline(138.5,41)(133.25,16.5)
\multiput(138.5,41)(-.033653846,-.157051282){156}{\line(0,-1){.157051282}}
%\end
%\emline(150.5,40.75)(153.25,17)
\multiput(150.5,40.75)(.03353659,-.28963415){82}{\line(0,-1){.28963415}}
%\end
%\emline(54.25,28.25)(54.5,25)
\multiput(54.25,28.25)(.03125,-.40625){8}{\line(0,-1){.40625}}
%\end
%\emline(56.5,28)(54,24.25)
\multiput(56.5,28)(-.03333333,-.05){75}{\line(0,-1){.05}}
%\end
%\emline(61.75,30)(62.5,27.75)
\multiput(61.75,30)(.0326087,-.0978261){23}{\line(0,-1){.0978261}}
%\end
%\emline(63,29.75)(62.75,27)
\multiput(63,29.75)(-.03125,-.34375){8}{\line(0,-1){.34375}}
%\end
%\emline(55.75,15.5)(57.75,16.5)
\multiput(55.75,15.5)(.0666667,.0333333){30}{\line(1,0){.0666667}}
%\end
%\emline(56.25,15.75)(57,15.25)
\multiput(56.25,15.75)(.05,-.0333333){15}{\line(1,0){.05}}
%\end
%\emline(57,15.25)(57.5,14.5)
\multiput(57,15.25)(.0333333,-.05){15}{\line(0,-1){.05}}
%\end
%\emline(20.25,15.5)(24,16.75)
\multiput(20.25,15.5)(.09868421,.03289474){38}{\line(1,0){.09868421}}
%\end
%\emline(21,15)(22.75,14.5)
\multiput(21,15)(.1166667,-.0333333){15}{\line(1,0){.1166667}}
%\end
\put(25,35.25){\line(1,0){2.25}}
%\emline(26,33.25)(28.25,35.5)
\multiput(26,33.25)(.03358209,.03358209){67}{\line(0,1){.03358209}}
%\end
%\emline(42.75,15.75)(43,16.75)
\multiput(42.75,15.75)(.03125,.125){8}{\line(0,1){.125}}
%\end
%\emline(42.75,15.5)(42.5,16.5)
\multiput(42.75,15.5)(-.03125,.125){8}{\line(0,1){.125}}
%\end
\put(134,25){\line(0,-1){2.5}}
%\emline(136,25)(134.75,23)
\multiput(136,25)(-.03289474,-.05263158){38}{\line(0,-1){.05263158}}
%\end
%\emline(151.5,25)(152.25,23.25)
\multiput(151.5,25)(.0326087,-.076087){23}{\line(0,-1){.076087}}
%\end
\put(153,25.75){\line(0,-1){2.5}}
%\emline(136.75,16.5)(139.75,17.5)
\multiput(136.75,16.5)(.1,.0333333){30}{\line(1,0){.1}}
%\end
%\emline(136.75,16)(139,15.25)
\multiput(136.75,16)(.0978261,-.0326087){23}{\line(1,0){.0978261}}
%\end
%\emline(146.75,42.5)(149.75,41.25)
\multiput(146.75,42.5)(.07894737,-.03289474){38}{\line(1,0){.07894737}}
%\end
%\emline(147,40.25)(150.25,40.75)
\multiput(147,40.25)(.2166667,.0333333){15}{\line(1,0){.2166667}}
%\end
\put(113.75,35.25){\line(1,0){3}}
%\emline(115,33.25)(117,35.25)
\multiput(115,33.25)(.03333333,.03333333){60}{\line(0,1){.03333333}}
%\end
\put(35.75,27.25){\makebox(0,0)[cc]{$\Delta$}}
\put(42.5,17.75){\makebox(0,0)[cc]{$o$}}
\put(26.25,18.75){\makebox(0,0)[cc]{$p$}}
\put(25.75,38.25){\makebox(0,0)[cc]{$q$}}
\put(51.25,22.75){\makebox(0,0)[cc]{$y^i$}}
\put(67.75,24){\makebox(0,0)[cc]{$y^{i+1}$}}
\put(57,10.25){\makebox(0,0)[cc]{$x_1^i$}}
\put(58,43.75){\makebox(0,0)[cc]{$x_2^i$}}
\put(52.75,38){\makebox(0,0)[cc]{$v_i$}}
\put(66,40.5){\makebox(0,0)[cc]{$v_{i+1}$}}
\put(57.5,19.5){\makebox(0,0)[cc]{$\Delta_i$}}
\put(115.5,40){\makebox(0,0)[cc]{$q$}}
\put(114.5,19.5){\makebox(0,0)[cc]{$p$}}
\put(121.25,28){\makebox(0,0)[cc]{$\Delta$}}
\put(131.5,30.25){\makebox(0,0)[cc]{$\bar y^k$}}
\put(157.75,27.5){\makebox(0,0)[cc]{$\bar y^{k+1}$}}
\put(141.25,11.5){\makebox(0,0)[cc]{$\bar x_1^k$}}
\put(141.5,44.5){\makebox(0,0)[cc]{$\bar x_2^k$}}
\put(143.25,23.75){\makebox(0,0)[cc]{$\Gamma_k$}}
\thinlines
%\emline(56.25,34.25)(62,39)
\multiput(56.25,34.25)(.040780142,.033687943){141}{\line(1,0){.040780142}}
%\end
%\emline(55.5,28.75)(62.25,35)
\multiput(55.5,28.75)(.036290323,.033602151){186}{\line(1,0){.036290323}}
%\end
%\emline(55.25,24.25)(61.75,30)
\multiput(55.25,24.25)(.038011696,.033625731){171}{\line(1,0){.038011696}}
%\end
%\emline(59,23.75)(62.75,27.75)
\multiput(59,23.75)(.033482143,.035714286){112}{\line(0,1){.035714286}}
%\end
%\emline(54.5,16.75)(56,18.25)
\multiput(54.5,16.75)(.03333333,.03333333){45}{\line(0,1){.03333333}}
%\end
%\emline(57,15.75)(62,21)
\multiput(57,15.75)(.033557047,.035234899){149}{\line(0,1){.035234899}}
%\end
%\emline(60.5,16.25)(63,18.5)
\multiput(60.5,16.25)(.03731343,.03358209){67}{\line(1,0){.03731343}}
%\end
%\emline(137.5,34.75)(146,40.75)
\multiput(137.5,34.75)(.047752809,.033707865){178}{\line(1,0){.047752809}}
%\end
%\emline(136.5,29.5)(150.5,40.5)
\multiput(136.5,29.5)(.0428134557,.0336391437){327}{\line(1,0){.0428134557}}
%\end
%\emline(138.75,27.5)(150.75,37)
\multiput(138.75,27.5)(.0425531915,.0336879433){282}{\line(1,0){.0425531915}}
%\end
%\emline(134.25,19.25)(138.25,22.5)
\multiput(134.25,19.25)(.041237113,.033505155){97}{\line(1,0){.041237113}}
%\end
%\emline(143.5,27.5)(151,33.25)
\multiput(143.5,27.5)(.043859649,.033625731){171}{\line(1,0){.043859649}}
%\end
%\emline(137.75,17.5)(137.25,16.5)
\multiput(137.75,17.5)(-.0333333,-.0666667){15}{\line(0,-1){.0666667}}
%\end
%\emline(135.5,16.5)(143.25,22.25)
\multiput(135.5,16.5)(.045321637,.033625731){171}{\line(1,0){.045321637}}
%\end
%\emline(147.75,26.75)(150.75,28.5)
\multiput(147.75,26.75)(.05769231,.03365385){52}{\line(1,0){.05769231}}
%\end
%\emline(142,17)(149.25,22.75)
\multiput(142,17)(.042397661,.033625731){171}{\line(1,0){.042397661}}
%\end
%\emline(146.5,17)(152.25,21)
\multiput(146.5,17)(.048319328,.033613445){119}{\line(1,0){.048319328}}
%\end
%\emline(151,17)(152.75,18.25)
\multiput(151,17)(.04605263,.03289474){38}{\line(1,0){.04605263}}
%\end
%\emline(138.25,39.25)(140.5,41)
\multiput(138.25,39.25)(.04326923,.03365385){52}{\line(1,0){.04326923}}
%\end
\put(35.5,6.25){\makebox(0,0)[cc]{Figure 7}}
\put(124,7.25){\makebox(0,0)[cc]{Figure 8}}
\end{picture}

\end{center}
\end{figure}

Consider the boxes $\Delta_i,\dots, \Delta_{j}$.  One can unite some neighbors to obtain larger
boxes $\Gamma_k$ ($k=1,\dots, m$ for some $m$) with contours $\bar x^k_1 (\bar y^k)^{-1} \bar x^k_2 \bar y^{k+1}$, where
$16|A|\le|\bar x^k_1| <32|A|$ (Fig. 8). Therefore $m\ge \frac{\theta n|A|-48|A|}{32|A|}> \frac{\theta n}{40}$.
We also derive $|\bar y^k|>\theta n|A|/2$ from the triangle inequality because the distance from $o$ to $(\bar y^k)_-$ is at least $\theta n|A|$.

Now by Lemma \ref{16}, we obtain $|\bar x^k_1|+|\bar x^k_2|>0.01 (\theta n|A|/2+\theta n|A| /2)= 0.01\theta n |A|$,
whence $|\bar x^k_2|>0.01\theta n|A| - 32|A|> 0.005\theta n |A|.$ It follows that $|q|\ge m(0.005\theta n|A|)>
\frac{\theta n}{40}(0.005\theta n|A|)>10^{-4}\theta^2|A|n^2\ge cn^2$ for $c= 10^{-4}\theta^2 $, and the lemma is proved.
\endproof

\bigskip

\section{Some non-simply-connected diagrams.}
We need a few types of non-simply-connected diagrams on surfaces.
They appear when one identifies some parts of the boundaries of simply-connected (= disc) diagrams. However if a simply-connected diagram $\Delta$ is not homeomorphic to the
standard disc, we can obtain undesired singularities after such  identifications. To avoid them, one can consider
only van Kampen diagrams homeomorphic to disc, where every cell is also homeomorphic to a disc. For this
goal, we suggest using so called $0$-cells
corresponding to trivial relations. By definition, every $0$-{\it edge} is labeled by the empty word $1$,
and we assign zero length to $0$-edges. The boundary contour of every $0$-{\it cell} (which is a disc too)
is of the form $e_1\dots e_n$, where all the edges $e_i$-s are $0$-edges or there are two edges $e_i$ and $e_j$
labeled by mutual inverse letters from the alphabet $\{a^{\pm 1}, b^{\pm 1}\}$, while the remaining
edges are $0$-edges.

Under such an agreement, one can use only surfaces without singularities. For instance,
every disc (spherical, annual, toric) diagram is now homeomorphic to  a disc (to a sphere, to an annulus, to a torus, resp.)
The van Kampen lemma (and other related lemmas) can be obviously reformulated for the diagrams with $0$-cells.
(See more details in Section 11 of \cite{O89}.) However instead of common arcs between $\cal R$-cells
$\Pi_1$ and $\Pi_2$, we consider now {\it contiguity subdiagrams} consisting of $0$-cells only and
having a contour of the form $p_1q_1p_2q_2$, where $|p_1|=|p_2|=0$, $q_1$ is a subpath of $\partial'\Pi_1$,
and $q_2$ is a subpath in $\partial'\Pi_2$ (Fig. 9). All the lemmas proved in the previous sections make
clear sense for diagrams with $0$-cells.

\begin{figure}
\begin{center}
% This is a LaTeX picture output by TeXCAD.
% File name: [emb5.pic].
% Version of TeXCAD: 4.3
% Reference / build: 30-Jun-2012 (rev. 105)
% For new versions, check: http://texcad.sf.net/
% Options on the following lines.
%\grade{\on}
%\emlines{\off}
%\epic{\off}
%\beziermacro{\on}
%\reduce{\on}
%\snapping{\off}
%\pvinsert{% Your \input, \def, etc. here}
%\quality{8.000}
%\graddiff{0.005}
%\snapasp{1}
%\zoom{4.0000}
\unitlength 1mm % = 2.845pt
\linethickness{0.4pt}
\ifx\plotpoint\undefined\newsavebox{\plotpoint}\fi % GNUPLOT compatibility
\begin{picture}(124,43.75)(10,0)
\thicklines
\put(23,30.25){\line(1,0){93}}
\put(23.25,27){\line(1,0){92}}
\put(23.75,24.25){\line(1,0){92}}
\put(23.5,30.25){\line(0,-1){5.75}}
%\emline(115.5,30.25)(115.25,24.25)
\multiput(115.5,30.25)(-.03125,-.75){8}{\line(0,-1){.75}}
%\end
\put(55,30.25){\line(0,-1){5.75}}
%\emline(86.25,30)(86.5,25)
\multiput(86.25,30)(.03125,-.625){8}{\line(0,-1){.625}}
%\end
\put(86.5,25.5){\line(0,-1){.75}}
%\emline(23.25,30.5)(16,39.75)
\multiput(23.25,30.5)(-.03372093,.043023256){215}{\line(0,1){.043023256}}
%\end
%\emline(115.25,30.5)(124,39.5)
\multiput(115.25,30.5)(.0336538462,.0346153846){260}{\line(0,1){.0346153846}}
%\end
%\emline(23.75,24.75)(15.5,16.75)
\multiput(23.75,24.75)(-.034663866,-.033613445){238}{\line(-1,0){.034663866}}
%\end
%\emline(115,24.5)(122.5,16.5)
\multiput(115,24.5)(.033632287,-.035874439){223}{\line(0,-1){.035874439}}
%\end
%\emline(74.25,32)(76.75,30.25)
\multiput(74.25,32)(.04807692,-.03365385){52}{\line(1,0){.04807692}}
%\end
%\emline(74,28.75)(77,30.5)
\multiput(74,28.75)(.05769231,.03365385){52}{\line(1,0){.05769231}}
%\end
%\emline(61.5,24.75)(63.75,25.5)
\multiput(61.5,24.75)(.0978261,.0326087){23}{\line(1,0){.0978261}}
%\end
%\emline(61.25,24.25)(63.5,23)
\multiput(61.25,24.25)(.05921053,-.03289474){38}{\line(1,0){.05921053}}
%\end
\put(18.25,27){\makebox(0,0)[cc]{$p_1$}}
\put(64.5,34.25){\makebox(0,0)[cc]{$q_1$}}
\put(119,26.75){\makebox(0,0)[cc]{$p_2$}}
\put(67.75,20.25){\makebox(0,0)[cc]{$q_2$}}
\put(40.25,37.75){\makebox(0,0)[cc]{$\Pi_1$}}
\put(94,18){\makebox(0,0)[cc]{$\Pi_2$}}
\thinlines
%\emline(26,42)(20,35.75)
\multiput(26,42)(-.033707865,-.03511236){178}{\line(0,-1){.03511236}}
%\end
\put(32.5,41.75){\line(-1,-1){9.75}}
%\emline(35.5,37.75)(28.75,31.5)
\multiput(35.5,37.75)(-.036290323,-.033602151){186}{\line(-1,0){.036290323}}
%\end
%\emline(39.25,35.25)(34.25,30.5)
\multiput(39.25,35.25)(-.035460993,-.033687943){141}{\line(-1,0){.035460993}}
%\end
\put(50.5,41){\line(-1,-1){10.25}}
%\emline(56.25,41.25)(46.5,31.75)
\multiput(56.25,41.25)(-.0345744681,-.0336879433){282}{\line(-1,0){.0345744681}}
%\end
%\emline(62.5,41)(51.75,30.75)
\multiput(62.5,41)(-.0353618421,-.0337171053){304}{\line(-1,0){.0353618421}}
%\end
%\emline(69,41.25)(65.75,38.5)
\multiput(69,41.25)(-.03963415,-.03353659){82}{\line(-1,0){.03963415}}
%\end
%\emline(60.75,33.5)(57.5,31)
\multiput(60.75,33.5)(-.04333333,-.03333333){75}{\line(-1,0){.04333333}}
%\end
%\emline(75.25,41.25)(69.75,36.25)
\multiput(75.25,41.25)(-.036912752,-.033557047){149}{\line(-1,0){.036912752}}
%\end
%\emline(81.75,42)(70.25,31)
\multiput(81.75,42)(-.0351681957,-.0336391437){327}{\line(-1,0){.0351681957}}
%\end
%\emline(88,42.25)(76.75,31.25)
\multiput(88,42.25)(-.0344036697,-.0336391437){327}{\line(-1,0){.0344036697}}
%\end
%\emline(95.25,43)(82.75,31.25)
\multiput(95.25,43)(-.0358166189,-.0336676218){349}{\line(-1,0){.0358166189}}
%\end
%\emline(102,43.75)(89.25,31.75)
\multiput(102,43.75)(-.0358146067,-.0337078652){356}{\line(-1,0){.0358146067}}
%\end
%\emline(108.75,43.25)(95.75,31.5)
\multiput(108.75,43.25)(-.0372492837,-.0336676218){349}{\line(-1,0){.0372492837}}
%\end
%\emline(114.25,43)(101.5,30.75)
\multiput(114.25,43)(-.0350274725,-.0336538462){364}{\line(-1,0){.0350274725}}
%\end
%\emline(119.75,42.5)(109.25,31.5)
\multiput(119.75,42.5)(-.0336538462,-.0352564103){312}{\line(0,-1){.0352564103}}
%\end
%\emline(29.5,23.5)(20.5,16.25)
\multiput(29.5,23.5)(-.041860465,-.03372093){215}{\line(-1,0){.041860465}}
%\end
%\emline(35.5,23.25)(26.25,16)
\multiput(35.5,23.25)(-.043023256,-.03372093){215}{\line(-1,0){.043023256}}
%\end
%\emline(42.5,24)(33,16.25)
\multiput(42.5,24)(-.041304348,-.033695652){230}{\line(-1,0){.041304348}}
%\end
%\emline(48.5,23.25)(39.75,16.5)
\multiput(48.5,23.25)(-.043532338,-.03358209){201}{\line(-1,0){.043532338}}
%\end
%\emline(55.5,23.25)(45.75,15.75)
\multiput(55.5,23.25)(-.043721973,-.033632287){223}{\line(-1,0){.043721973}}
%\end
%\emline(61.25,23.5)(51.5,15.5)
\multiput(61.25,23.5)(-.040966387,-.033613445){238}{\line(-1,0){.040966387}}
%\end
%\emline(63,20.5)(57.75,16)
\multiput(63,20.5)(-.039179104,-.03358209){134}{\line(-1,0){.039179104}}
%\end
%\emline(66.25,18.25)(63.25,15.5)
\multiput(66.25,18.25)(-.03658537,-.03353659){82}{\line(-1,0){.03658537}}
%\end
%\emline(77.25,23.5)(68.75,15.75)
\multiput(77.25,23.5)(-.036956522,-.033695652){230}{\line(-1,0){.036956522}}
%\end
%\emline(83.25,24)(74.25,16)
\multiput(83.25,24)(-.037815126,-.033613445){238}{\line(-1,0){.037815126}}
%\end
%\emline(89.25,23.25)(79.5,15.75)
\multiput(89.25,23.25)(-.043721973,-.033632287){223}{\line(-1,0){.043721973}}
%\end
%\emline(88.75,17.5)(85.25,14.75)
\multiput(88.75,17.5)(-.04268293,-.03353659){82}{\line(-1,0){.04268293}}
%\end
%\emline(103.5,24)(98.75,21)
\multiput(103.5,24)(-.05337079,-.03370787){89}{\line(-1,0){.05337079}}
%\end
%\emline(93.5,15.5)(91.25,13.75)
\multiput(93.5,15.5)(-.04326923,-.03365385){52}{\line(-1,0){.04326923}}
%\end
%\emline(108.5,23.5)(97.75,14.75)
\multiput(108.5,23.5)(-.0413461538,-.0336538462){260}{\line(-1,0){.0413461538}}
%\end
%\emline(114.5,23.75)(104,14.5)
\multiput(114.5,23.75)(-.0381818182,-.0336363636){275}{\line(-1,0){.0381818182}}
%\end
%\emline(117.75,20.75)(110.75,15)
\multiput(117.75,20.75)(-.040935673,-.033625731){171}{\line(-1,0){.040935673}}
%\end
%\emline(120.5,18.25)(118.5,15.75)
\multiput(120.5,18.25)(-.03333333,-.04166667){60}{\line(0,-1){.04166667}}
%\end
\put(63.5,7){\makebox(0,0)[cc]{Figure 9}}
\end{picture}

\end{center}
\end{figure}

If a diagram $\Delta$ has a simple path $x$ connecting two vertices $o_1$ and $o_2$, and a diagram $\Delta'$
results from $\Delta$ after some amalgamation or cancellation of cells, which does not affect $o_1$ and $o_2$,
then using $0$-cells one can construct a simple path $x'$ connecting the same vertices in $\Delta$
with the label $Lab(x')$ freely equal to $Lab(x)$. (See details in Section 13.5 of \cite{O89}.)

One more distinction in comparison with simply-connected case is that a cell can  be compatible
with itself in a non-simply-connected diagram $\Delta$. This means that in the definition of compatibility,
we now allow the equality $\Pi_1=\Pi_2$,
i.e., a cell $\Pi_1$ is compatible with itself if the path $o_1-o_2$
together with a subpath of $\partial\Pi$ connecting the entire vertices $o_1$ and $o_2$, give a closed
path, which is not $0$-homotopic in $\Delta$. This path is labeled by a $G$-word.

We call a diagram $\Delta$ {\it singular} if it contains a simple closed path, which
is not $0$-homotopic in $\Delta$ but has trivial in the free group label.

\begin{lemma} \label{holes} Let $\Delta$ be a reduced diagram on a sphere with $n\ge 2$ holes, i.e.,
$\Delta$ has $n$ boundary components $p_1, \dots, p_n$ with clockwise labels $P_1,\dots, P_n$.

(1) If $\Delta$ is singular, then there are a positive integer  $m\in [1,n-1]$,   indices $i_1<\dots<i_m$,
and words $Q_1,\dots, Q_m$ such that $\prod_{j=1}^m Q_jP_{i_j}Q_j^{-1}=1$ in $H$.

(2) If $\Delta$ is non-singular and the words $P_1,\dots,P_n$ are $G$-words, then
the label of arbitrary path connecting in $\Delta$ any two entire vertices of the cells
or of the boundary components, is freely equal to a $G$-word.
\end{lemma}

\proof (1) We have a simple closed path $q$, such that $Lab(q)=1$ in $F(a,b)$, and $q$ bounds
a subdiagram $\Gamma$ with $m$ holes, where $1\le m<n$. The holes are bounded by some $p_{i_1},\dots, p_{i_m}$.
Then, after cancellations in $q$, $\Gamma$ becomes a diagram on a sphere with $m$ holes.
The standard application of the version of van Kampen -- Schupp Lemma (see Chapter V in \cite{LS})
provides us with an equality $\prod_{j=1}^m Q_jP_{i_j}Q_j^{-1}=1$ in $H$.

(2) Let us consider the bigger set of relations $\cal Q$ which consists of all $G$-words. Then we obtain a spherical
$\cal Q$-diagram $\Gamma$ from $\Delta$ if we patch up every hole by a new cell. It suffices to prove
the property (2) for $\Gamma$. We will distinguish the old cells, i.e. the cells of $\Delta$ and the new ones. By Lemma \ref{Grin} (1) (where the length of the boundary is $0$), the
diagram $\Gamma$ is not reduced. Hence it has two different compatible cells $\Pi_1$ and $\Pi_2$
which can be replaced by a single $\cal Q$-cell $\Pi$ in the spherical diagram $\Gamma'$ with fewer cells.

Note that both $\Pi_1$ and $\Pi_2$ cannot be old since then the original diagram could not be reduced.
If one of them is old and another one is new, then we say that $\Pi$ is new. This cell $\Pi$ cannot have trivial
label in the free group since the original diagram $\Delta$ was non-singular.

If both $\Pi_1$ and $\Pi_2$
are new, we say that $\Pi$ is also new but with multiplicity $2$. Again $\Pi$ is labeled by a freely non-trivial $G$-word unless
$n=2$, because $\Delta$ was non-singular. If $\Pi$ is non-trivial, i.e., $\Pi_1$ and $\Pi_2$ do not annihilate, then $\Pi$ will have at least one entire vertex.

If this cell $\Pi$ is non-trivial, then it suffices to prove the property (2)  for the
diagram $\Gamma'$ with fewer cells because the cells $\Pi_1$ and $\Pi_2$ of $\Delta$ were compatible
and the compliment of a subpath between two entire vertices in $\partial\Pi$ is also labeled by a $G$-word.

Proceeding this way, we obtain a chain $\Gamma = \Gamma^{(0)}\to\Gamma'\to\dots\to \Gamma^{(s)}$, where $\Gamma^{(s)}$ consists
of new and old $\cal Q$-cells. Every step $\Gamma^{(i-1)}\to\Gamma^{(i)}$ decreases the number of old cells or
replaces two new cells of multiplicities $m_1$ and $m_2$ by a single cell of multiplicity $m_1+m_2<n$.
Thus, soon or later one achieves a transition $\Gamma^{(i-1)}\to\Gamma^{(i)}$, where two new cells $\Pi_1$ and $\Pi_2$ of
$\Gamma^{(i-1)}$ with multiplicities $m_1$ and $m_2$ are compatible and $m_1+m_2=n$.

If $\Pi_1$ and $\Pi_2$ annihilates in $\Gamma^{(i)}$, then $\Gamma^{(i)}$ has no old cells by Lemma \ref{Grin} (1), being a spherical
reduced diagram. The same is true for $\Gamma^{(i-1)}$, and so $\Gamma^{(i-1)}$ has no entire vertices except for
the vertices of $\Pi_1$ and $\Pi_2$. Since these two cells are compatible, the hypothesis (b) holds for $\Gamma^{(i-1)}$.
Hence it holds for $\Delta$.

If $\Pi_1$ and $\Pi_2$  nontrivially amalgamate in $\Gamma^{(i)}$, then one can remove their amalgamation
$\Pi$ from $\Gamma^{(i)}$ and obtain a disc diagram $E$  over $H$, whose contour is an $G$-word.
Again one can remove the cells from $E$ one-by-one changing the contour of $E$ and stop when a diagram
$E^{(i)}$ has only one (Grindlinger) cell.

\endproof

\begin{lemma} \label{tor} Let $\Delta$ be a reduced toric diagram.

(1) If $\Delta$ is singular, then the boundary labels of all closed
paths starting at a fixed vertex of $\Delta$ belong to the same  cyclic subgroup of $H$.

(2) If $\Delta$ is not singular, then the label of arbitrary path connecting the entire
vertices of its cells is freely equal to a $G$-word.

\end{lemma}

\proof (1) There is a simple closed path $x$ in $\Delta$
such that it represents non-trivial element of the fundamental group $G_T$ of the torus,
but $Lab(x)$ is freely equal to $1$. We can find a simple closed path $y$, such that
$x$ and $y$ represent a canonical generating set of $G_T.$  In this case the
subgroup $\langle Lab(x), Lab(y)\rangle $ of $G$ is cyclic and so the image of $G_T$ is cyclic
under the homomorphism $Lab$ from $G_T$ to $H$. The first statement is
proved since changing the origin, one replaces this subgroup by a conjugate one.

(2) Let us contract every $0$-edge of $\Delta$ to a vertex and, thus, contract every $0$-cell to a vertex or
to an edge. We obtain a reduced diagram $\Gamma$ on a surface homeomorphic to a torus since no closed path, which is not $0$-homotopic,
will be contracted to a vertex. So it suffices to prove the statement for $\Gamma$.

If for every cells $\Pi_1$ and $\Pi_2$, the length of any
common arc $q$ of their reduced  boundaries
$\le 0.05|\partial'\Pi_1|$,
then every cell of $\Gamma$ is an $n$-gon with $n>20$, but the torus cannot be
tessellated by such polygons, as this follows from Euler's formula, a contradiction.
Otherwise one can find a pair of
compatible cells $\Pi_1$ and $\Pi_2.$ In fact $\Pi_1=\Pi_2$ since the diagram $\Gamma$ is reduced.
In other words, we have a self-compatible cell $\Pi$ in $\Gamma$.

A part of $\partial\Pi$ and the path $x$ defining the self-compatibility, form a closed path $p$, which is
not contractible to a vertex along $\Gamma.$ Therefore, if we cut the torus along $p$, we obtain an annular diagram
$E$ whose boundary labels are freely equal to $G$-words. The diagram $E$ is not singular
since $\Gamma$ is non-singular. The statement (2) holds for $\Gamma$ iff it holds
for $E$ since $\Pi$ is a self-compatible cell. However the required property of $E$ has
been already proved in Lemma \ref{holes} (2). This completes the proof of Lemma \ref{tor}.
\endproof

\section{Algebraic properties of the embedding $\gamma$.}

\begin{lemma} \label{comm} If a word $W$ is free and $WV=VW$ in $H$ for some word $V$, then the subgroup
$\langle W,V \rangle $ of $H$ is cyclic.
\end{lemma}

\proof By van Kampen's Lemma, we have a disc diagram $\Delta$ over $G$ with
contour $p_1q_1p_2q_2$, where $Lab(p_1)\equiv Lab(p_2)^{-1}\equiv W$ and
$Lab(q_1)\equiv Lab(q_2)^{-1}\equiv V.$ Identifying $p_1$ with $p_2^{-1}$ and $q_1$
with $q_2^{-1}$, we have a toric diagram $\Gamma.$ Let $\Gamma_0$ be a reduced
toric diagram resulting from $\Gamma$ after possible amalgamations and annihilations
of some cells. Using auxiliary $0$-cells, one may assume that it also has closed paths with labels freely equal to $V$ and $W$.
Consider now two cases.

If $\Gamma_0$ is singular, then
the statement follows from Lemma \ref{tor} (1).
If $\Gamma_0$ is nonsingular, then by Lemma \ref{tor} (2), the words $V$ and $W$
are freely conjugate to some $G$-words, a contradiction.

\endproof

\begin{lemma}\label{notinG} If $A$ is a free word, then so is $A^t$ for $t\ne 0$.
\end{lemma}
\proof One may assume that $A$ is minimal and $t\ge 2$. Arguing by contradiction,
assume that the word $A^t$ is not free, and therefore there is a reduced annular
diagram $\Delta$ with contours $p_1$ and $p_2$ labeled by $A^t$ and by a $G$-word,
respectively. This diagram is non-singular since $A^t\ne 1$ in $H$ by Lemma \ref{quasi} (1).
So, as in the proof of Lemma \ref{tor}, one may assume that $\Delta$ has no $0$-edges.
If $\Delta$ has no cells, then we can compare the labels of $p_1$ and $p_2$ using
Lemma \ref{0.55} (1), and obtain a contradiction since $t\ge 2>1.1.$ Therefore $\Delta$
has an $\cal R$-cell $\Pi$.

A cell $\Pi$ cannot have a common boundary arc of length $\ge 0.05|\partial'\Pi|$ with $p_2$ since this would make possible
to modify $p_2$ using Lemma \ref{dia} (2) and to reduce the number of cells in $\Delta$. By Lemma \ref{dia} (1),
$\Pi$ has no common boundary arcs of length $\ge 0.05|\partial'\Pi|$ with other cells. The unique maximal common boundary
arc of $\Pi$ and $p_1$ has length $\le 0.55|\partial'\Pi|$ by Lemmas \ref{1touch} and \ref{0.55}. Thus, either the polygon  $\Pi$ has at least $11$ sides,
or the cell $\Pi$ is self-compatible. The former case for all the cells of $\Delta$ gives a contradiction with
Euler' formula. If we have the latter case for some $\Pi$, then we cut up $\Delta$
by a path labeled by a $G$-word, as we did that in the proof of Lemma \ref{tor} (2). This gives an annular subdiagram
$\Delta_1$, where one contour is $p_1$, and another one is again labeled by an $G$-word. Since $\Delta_1$
has fewer $\cal R$-cells than $\Delta$, the statement is proved by induction.

\endproof

\begin{lemma} \label{st} Let $A$ be a free word. Then any equality $XA^sX^{-1}=A^t$ implies that
$|s|=|t|$.
\end{lemma}

\proof One may assume that the word $A$ is minimal. Since $X^kA^{s^k}X^{-k}=A^{t^k}$ for any natural $k$,
 we have $|s^k|>0.2|t^k|$ for every $k$ by Lemma \ref{quasi} (2), whence $|s|\ge |t|.$ By the symmetry, we also have $|t|\ge |s|$. The lemma is proved.
 \endproof

If a path $q$ has an $A$-periodic label, then one can select vertices along $q$
 such that they divide $q$ into subpaths of length $\le |A|$ and any two
 vertices of this system are connected by a subpath of $q$ labeled by a power of $A$.
 The vertices of this system will be called the {\it phase} vertices of $q$.

\begin{lemma} \label{P1is} Let $A$ be a minimal word. Then for every constant $C$,
there is a natural number $m_0=m_0(C, |A|)$ with the following property. If $P_1Q_1P_2Q_2=1$
in $H$, where $|P_1|, |P_2|\le C$,  $Q_1$ and $Q_2^{-1}$ are $A$- or $A^{-1}$-periodic words
starting with $A^{\pm 1}$ and $|Q_1|\ge m_0|A|$, then $P_1A^tP_1^{-1} =A^{\pm t}$ in $H$
for some $t\ne 0$.
\end{lemma}

\proof Consider a reduced diagram with contour $p_1q_1p_2q_2$, where $Lab(p_i)\equiv P_i$, $Lab(q_i)\equiv Q_i$
($i=1,2$). If $m_0$ is large enough, then by Lemma \ref{polosa} (1), there is a phase vertex $o$ on $q_1$
connected with a phase vertex $o'$ on $q_2$ by
a path of length $<2800|A|+2|A|<3000|A|$. Moreover by Lemma \ref{polosa} (2), we may assume that the number of such phase
vertices is sufficiently large to guarantee that there are different phase vertices  $o(1)$ and $o(2)$  of $q_1$  connected with some phase vertices of $q_2$ by paths having the same label $Z$. Therefore there is a closed path in $\Delta$
with label $Z^{-1}A^tZ A^s,$ where $t\ne 0$. So this product is trivial in $H$ and $s=\pm t$ by Lemma \ref{st}.
Since the vertices $(p_1)_{\pm}$ are phase vertices two, we also have $Z^{-1}=A^kP_1A^l$ in $H$ for some
integers $k$ and $l$. Therefore $(A^kP_1A^l)A^t(A^kP_1A^l)^{-1}A^{\pm t}=1$ in $H$, which proves the
lemma.
\endproof

For a free element $A$, we introduce its {\it elementary closure} $E(A)$:
$$ E(A) = \{X\in H\mid  XA^tX^{-1}=A^{\pm t}\;\; for \;\;some \;\; t=t(X)\ne 0\}$$
Clearly, $E(A)$ is a subgroup of $H$ containing the cyclic subgroup $\langle A \rangle.$

\begin{lemma} \label{fin} The index $[E(A):\langle A\rangle]$ is finite.
\end{lemma}
\proof  Since $E(BAB^{-1})=BE(A)B^{-1}$ for any $B$, we may assume that $A$ is a minimal word.
Let $X\in E(A)$. Since $XA^tX^{-1}=A^{\pm t}$ for some $t>0$, we have
$XA^mX^{-1}=A^{\pm m}$ for any large $m$, which is a multiple of $t$. We consider a reduced
disc diagram $\Delta$ for the equality $XA^{2m}X^{-1}=A^{\pm 2m}$ with contour
$p_1q_1p_2q_2$, where $Lab(p_1)\equiv Lab(p_2)^{-1}\equiv X$, $Lab(q_1)\equiv A^{2m}$
and $Lab(q_2)\equiv A^{\mp 2m}$.

Let $o$ be the middle phase vertex of $q_1$ corresponding to the factorization $A^mA^m$ of $Lab (q_1)$.
There is a path connecting $o$ with the middle phase vertex $o'\in q_2$ having label $A^{-t}X^{-1}A^{\mp t}=X^{-1}$
in $H$. On the other hand, if $m=m(|X|, |A|)$ is very large, then by Lemma \ref{polosa} (2), the vertex
$o$ can be connected with a phase vertex $o''$ on $q_2$ by a path having label $Y$ of length
$<3000|A|$. So the closed path $o-o'-o''-o$ gives the equality $X^{-1}A^kY^{-1}=1$ in $H$, where
the section $o'-o''$ of $q_2$ is labeled by $A^k$.

Hence $X=A^kY^{-1}$ in $H$, and so every element $X\in E(A)$ belongs to a right coset $\langle A\rangle Y^{-1}$
with $|Y|<3000|A|$. Therefore the set of right cosets of  $\langle A\rangle $ in $E(A)$ is finite,
as required.
\endproof

We define $ E^+(A) = \{X\in G\mid  XA^tX^{-1}=A^{t}\;\; for \;\;some \;\; t\ne 0\}$. It is easy to see
that $E^+(A)$ is a subgroup of $E(A)$ of index $\le 2$.

\begin{lemma} \label{cycl} The subgroup $E^+(A)$ is infinite cyclic for a free word $A$.
\end{lemma}

\proof Let $X\in E^+(A)\backslash\{1\}$. Then for some $t>0$, the subgroup
$\langle A^t,X\rangle$ is infinite cyclic by Lemmas \ref{notinG} and \ref{comm}. Therefore $X$ has infinite order.
Thus, the group $E^+(A)$ is torsion free and contains the cyclic group $\langle A\rangle$
of finite index by Lemma \ref{fin}. It follows that $E^+(A)$ is cyclic itself  by Schur's theorem (See \cite{R}, 10.1.4.)
\endproof

\begin{lemma} \label{free} Let $A$ be a minimal word and $B\in H\backslash E(A)$.
Then for sufficiently large $m$ the powers $A^m$ and $BA^mB^{-1}$ generate
a free subgroup of rank $2$ in $H$.
\end{lemma}

\proof Proving by contradiction, we have a nontrivial relation $r(A^m, BA^mB^{-1})=1,$
which gives a reduced diagram $\Delta$ with contour $p_1q_1\dots p_nq_n$, where
$Lab(p_i)\equiv B^{\pm 1}$ and $Lab(q_i)\equiv A^{k_im}$ for some integer $k_i\ne0$
($i=1,\dots,n$).

If $m$ is large enough, we can apply Lemma \ref{pnqn} to $\Delta$ and obtain a subdiagram
$\Gamma$ with boundary label of the form $P_1Q_1P_2Q_2$, where $Lab(P_1)\equiv B^{\pm 1},$
$|P_2|<15|A|$, $Q_1$ and $Q_2$ are periodic words with periods $A^{\pm 1}$, $|Q_1|>m|A|/3$ and
the words $Q_1$ and $Q_2^{-1}$ have prefixes $A$ or $A^{-1}$.

Then we may apply Lemma \ref{P1is} to the equality $P_1Q_1P_2Q_2=1$ in $H$ if $m/4$ is greater
than the constant $m_0$ depending on the lengths of the words $A$ and $B$. It says that
$B\in E(A)$, a contradiction.
\endproof

\bigskip

\section{Proofs of the theorems.}

{\bf Proof of Theorem \ref{t1}.} Assume that $K$ is not conjugate to a subgroup of $G.$
At first we also assume that $K$ has a free element $A$. We may replace $K$ by a conjugate
subgroup and assume that $A$ is minimal. If $K$ has an element $B\notin E(A)$,
then $K$ has a free subgroup of rank $2$ by Lemma \ref{free}.

Suppose $K\le E(A)$. The subgroup $E(A)$ has an infinite cyclic subgroup $E^+(A)$ of
index $\le 2$ by Lemma \ref{cycl}. But $E(A)$ has no elements of order $2$ in its center
by Lemma \ref{comm} applied to $A$. Hence $E(A)$ is either an infinite cyclic or an infinite
dihedral group. So is its subgroup $K$.

If $K$ is an infinite dihedral group, then it is generated by two involutions conjugated to
some elements from $H$, since free elements have infinite order by Lemma \ref{quasi} (1).

If $K$ is infinite cyclic, then $K\cap hGh^{-1}=\{1\}$ by Lemma \ref{notinG}.

Assume now that the subgroup $K$ has no free elements. Again, passing to a conjugate
subgroup, we may assume that $K$ contains a nontrivial element $A$ from $G$.
Arbitrary nontrivial element from $K\backslash\langle A\rangle$ is equal in $H$ to $CBC^{-1}$, where $B\in G\backslash\{1\}$.
Since $ACBC^{-1}\in K$, this product is conjugate to a non-trivial element $D$ from $G$.

Thus we have an annular diagram $\Delta$ for the conjugation of $ACBC^{-1}$ and $D$. Identifying
the boundary subpaths labeled by $C$ and $C^{-1}$, we obtain a diagram $\Delta_1$ on a sphere
with three holes bounded by closed paths $p_1$, $p_2$ and $p_3$ with labels $A, B$ and $D$, resp.
The vertices $(p_1)_-$ and $(p_2)_-$ are connected by a simple path $x$ labeled by $C$.

Let $\Delta_2$ is a reduced form of $\Delta_1$. It is non-singular since the words $A, B$ and $D$
are non-trivial in $H$. The vertices $(p_1)_-$ and $(p_2)_-$ can be connected by a simple path
$y$ such that $Lab(y)=Lab(x)=C$ in $H$. It follows now from Lemma
\ref{holes} (2), that $C$ is an $G$-word. Hence $CBC^{-1}\in G$, and so $K\le G.$
$\Box$
\medskip

{\bf Proof of Theorem \ref{t2}.} The statement follows from Lemma \ref{quasi} (1)
because any free element of $H$ is conjugate to a minimal element. $\Box$
\medskip

{\bf Proof of Theorem \ref{t4}.} Since a conjugation by a fixed element is a bi-Lipschitz mapping $H\to H$,
we may assume that the element $g$ is presented by a minimal word $A$. Now consider a path
$p$ from $g^n$ to $g^{-n}$ labeled by $A^{-2n}$. One can construct a reduced diagram $\Delta$
corresponding to the closed path $pq$. Since the distance in the Cayley graph between any two
vertices of $pq$ do not exceed the distance between the corresponding vertices of $\partial\Delta$,
we draw the required inequality from Lemma \ref{sq}. $\Box$

\medskip

{\bf Proof of Theorem \ref{t3}.}
We will prove both parts of the theorem using simultaneous induction on $n$ with obvious
base $n=1.$ Let us start with the second claim.

(2) Assume that we have a reluctant family $g_1,\dots,g_n\in G$ and the elements $x_1,\dots, x_n\in H$ such
that $\prod_{i=1}^n x_ig_ix_i^{-1}=1$.
We can construct a diagram $\Delta$ with boundary label $\prod_{i=1}^n X_iP_iX_i^{-1}$, where
the reduced words $P_1,\dots, X_1,\dots$ represent the elements $g_1,\dots,x_1,\dots$. Then for every $i$,
we identify the subpaths of $\partial\Delta$ labeled by $X_i$ and $X_i^{-1}$ and obtain a
diagram $\Delta_1$ on a sphere with $n$ holes. The $n$ contours $p_1,\dots p_n$ of $\Delta_1$
are labeled by $P_1,\dots P_n$. Now the path $x_i$ connects some vertex $o$ with the vertex $(p_i)_-$.
If one moves $o$ to another point, say to $(p_1)_-$, then all the labels $X_i$ are multiplied from the left
by the same word. Hence it suffices to prove that the label of every path $z_i$ connecting $(p_1)_-$
with $(p_i)_-$ is freely equal to a $G$-word.

Denote by $\Delta_2$ the reduced form of $\Delta_1.$ The vertex $(p_1)_-$ can be connected
in $\Delta_2$ with each $(p_i)_-$ by a path $u_i$ such that $Lab(u_i)=Lab(z_i)$ in $H$.

Note that the diagram $\Delta_2$ is non-singular. Indeed, otherwise we should have a simple closed path
 $q$ which bounds a subdiagram with $m$ holes, where $1\le m\le n-1$, and $Lab(q)=1$ in the free group.
But this means that for some $1\le i_1\le\dots\le i_m\le n$, the product of conjugacy classes $\prod_{j=1}^m g_{i_j}^H$ contains the identity element. Then applying the statement (1) for $m<n$, we conclude that the product $\prod_{j=1}^m g_{i_j}^G$
also contains $1$, which is impossible since the family $(g_1,\dots, g_n)$ is reluctant.

 Since the paths $u_i$ connect
some pairs of entire vertices in $\Delta_2$, $Lab(u_i)\in G$ by Lemma \ref{holes} (2), as desired.

(1) Assuming that $\prod_{i=1}^n x_ig_ix_i^{-1}=1$, we construct the diagrams $\Delta$, $\Delta_1$ and $\Delta_2$
as in the proof of the part (2).

If $\Delta_2$ is non-singular, then as in the proof of the claim (2), we obtain that the elements
$x_1,\dots, x_n$ belong to a coset $hG$ for some $h\in H$. Therefore $h^{-1}x_1,\dots, h^{-1}x_n\in G$
and  $$\prod_{i=1}^n (h^{-1}x_i)g_i(h^{-1}x_i)^{-1}=h^{-1}(\prod_{i=1}^n x_ig_ix_i^{-1})h=1$$.

If $\Delta_2$ is singular, then we have a simple closed path $q$ with trivial in the free group label,
and $q$ cut up $\Delta_2$. One of the two obtained parts gives us an equality of the form $\prod_{j=1}^m y_j g_{i_j}y^{-1}_j=1$
in $H$ for some $m\in [1,n-1]$ and $y_j\in H$. Another one gives an equality $\prod_{j=1}^{n-m} z_j g_{k_j}z^{-1}_j=1$, where
$\{k_1,\dots,k_{n-m}\}=\{1,\dots,n\}\backslash\{i_1,\dots, i_m\}$. In other words, both products
$\prod_{j=1}^m g_{i_j}^H$ and $\prod_{j=1}^{n-m} g_{k_j}^H$ contain $1$. By the statement (1), the
same is true for the products $\prod_{j=1}^m g_{i_j}^G$ and $\prod_{j=1}^{n-m} g_{k_j}^G.$ Hence
the product of $n$ conjugacy classes $g_{i}^G$  of the group $G$ contains the identity too, as required.
$\Box$

\bigskip

{\bf Acknowledgements.} The author is thankful to Mark Sapir who has given an impetus for
the present research.

\bigskip

\bigskip

{\bf Alexander A. Olshanskii:} Department of Mathematics, Vanderbilt University, Nashville
37240, U.S.A., and Moscow State University, Moscow 119991, Russia.

{\it E-mail}: alexander.olshanskiy@vanderbilt.edu

\end{document}